\documentclass[11pt]{amsart}

\usepackage{amssymb,amsfonts,latexsym}
\usepackage[centertags]{amsmath}
\usepackage{amsfonts}
\usepackage{amssymb}
\usepackage{amsthm}
\usepackage[all]{xy}
\usepackage{graphicx,color}

\newtheorem{cont}{cont}[section]
\newtheorem{theorem}[cont]{Theorem}
\newtheorem{proposition}[cont]{Proposition}
\newtheorem{lemma}[cont]{Lemma}
\newtheorem{corollary}[cont]{Corollary}
\newtheorem{definition}[cont]{Definition}

\newtheorem*{Enunciato*}{Enunciato}
\newtheorem{conj}[cont]{Conjecture}

\numberwithin{equation}{section}
\newtheorem{remark}[cont]{Remark}
\newtheorem{example}[cont]{Example}
\newtheorem*{not*}{Notation}

\newcommand{\cA}{\mathcal{A}}

\newcommand{\cC}{\mathcal{C}}
\newcommand{\cD}{\mathcal{D}}

\newcommand{\shE}{\mathcal{E}}
\newcommand{\PP}{\mathbb{P}}
\newcommand{\LL}{\mathbb{L}}
\newcommand{\ZZ}{\mathbb{Z}}

\newcommand{\QQ}{\mathbb{Q}}

\newcommand{\odi}[1]{\mathcal{O}_{#1}}

\DeclareMathOperator{\Hl}{H} 
 
\DeclareMathOperator{\rk}{rk} \DeclareMathOperator{\Hom}{Hom}

\DeclareMathOperator{\depth}{depth}

\DeclareMathOperator{\Proj}{Proj} 
\DeclareMathOperator{\codim}{codim}
\DeclareMathOperator{\coker}{coker}
\DeclareMathOperator{\Ext}{Ext} \DeclareMathOperator{\pd}{pd}
\DeclareMathOperator{\Tor}{Tor} 
 \DeclareMathOperator{\im}{im}
\DeclareMathOperator{\id}{id}

\begin{document}
\title{Schur powers of the cokernel of a graded morphism}

\author[Jan O.\ Kleppe, Rosa M.\ Mir\'o-Roig]{Jan O.\ kleppe, Rosa M.\
Mir\'o-Roig$^{*}$}
\address{OsloMet-Oslo Metropolitan University,
         Faculty of Technology, Art and Design,
         PB. 4, St. Olavs Plass, N-0130 Oslo,
         Norway}
\email{JanOddvar.Kleppe@oslomet.no}
\address{Universitat de Barcelona,
Departament de Matem\`{a}tiques i Inform\`{a}tica, Gran Via de les Corts Catalanes
585, 08007 Barcelona, SPAIN } \email{miro@ub.edu}

\date{\today}
\thanks{$^*$ Partially supported by MTM2016-81735-P and PID2019-104844GB-I00}

\begin{abstract} Let $\varphi: F\longrightarrow G$ be a graded morphism
  between free $R$-modules of rank $t$ and $t+c-1$, respectively, and let
  $I_j(\varphi)$ be the ideal generated by the $j \times j$ minors of a matrix
  representing $\varphi$. In this paper: (1) We show that the canonical
  module of $R/I_j(\varphi)$ is up to twist equal to a suitable Schur power
  $\Sigma ^I M$ of $M=\coker (\varphi ^*)$; thus equal to $\wedge ^{t+1-j}M$
  if $c=2$ in which case we find a minimal free $R$-resolution of
  $\wedge ^{t+1-j}M$ for any $j$, (2) For $c = 3$, we construct a free
  $R$-resolution of $\wedge ^2M$ which starts almost minimally (i.e. the first
  three terms are minimal up to a precise summand), and (3) For $c \ge 4$, we
  construct under a certain depth condition the first three terms of a free
  $R$-resolution of $\wedge ^2M$ which are minimal up to a precise summand. As
  a byproduct we answer the first open case of a question posed by Buchsbaum
  and Eisenbud in \cite[pg. 299]{BE75}. \end{abstract}

\maketitle


\section{Introduction}
The study of Schur functors and determinantal ideals has a long history that
goes back to the seminal work of Lascoux \cite{L} who used Schur functors to
construct (in characteristic 0) a free graded resolution of determinantal
ideals. More precisely, let $\varphi: F\longrightarrow G$ be a graded
$R$-morphism between free $R$-modules of rank $t$ and $t+c-1$, respectively.
Set $M=\coker (\varphi ^*)$. It is an interesting open problem in commutative
algebra to determine a graded free minimal
$R$-resolution of the cokernel
of the natural map,
$$\varphi _{r,s}^{*}:\bigwedge ^rG^*\otimes \bigwedge ^sF^*\longrightarrow
\bigwedge ^{r+s}F^*, \quad f\otimes g \mapsto \bigwedge ^r\varphi ^*(f)\wedge
g,$$
and in particular to find a graded free minimal $R$-resolution of the
$p$-exterior power, $\bigwedge^p M:= \coker (\varphi^*_{1,p-1})$, of $M$, as
well as to describe a graded free minimal $R$-resolution of any Schur power
$\Sigma ^I M$. So far only few particular cases have been described, and the
most relevant contribution to this problem is due to Lascoux who using Schur
complexes constructed a minimal free resolution of
$\coker(\varphi ^*_{r,t-r})=R/I_r(\varphi )$. In \cite{BE75},
Buchsbaum and Eisenbud solves the case  $\coker (\varphi^*_{p,0})$ and in
page 299 they predict the first terms of a graded free minimal
$R$-resolution of $\coker (\varphi^*_{p-1,1})$. No contributions have been
achieved until now and we devote this paper to study the first open cases and solve the case $p=c=2$  as well as the case $p=2$ and $c=3$ up to an explicitly described summand which seems to split off. Also the case $p=2$, $c\ge 4$ is solved under a certain depth condition, up to a concrete summand. More precisely, this work has two goals: (1) To
show that the canonical module of any determinantal ring $R/I_{i}(\varphi )$ is up to twist equal to a
suitable Schur power of the cokernel $M$ of $\varphi ^*$, thus equal to $\wedge ^{t+1-i}M$ if $c=2$; and (2) To construct
a free resolution of $\wedge ^2 M$ for the cases $c=2,3$ and, hence, solve Buchsbaum-Eisenbud question for $p=c=2$, and for $p=2$ and $c=3$ up to a concrete summand. Moreover, for $c=2$, we find a minimal $R$-free resolution of $\wedge ^p M$ for any $p$.

As a main tool we use the resolution of determinantal ideals $I_{i}(\varphi )$
constructed by Lascoux \cite{L} using Schur functors. For sake of completeness
we quickly review the construction. As a byproduct we obtain the minimal
resolution of the normal module of a standard determinantal ideal of
codimension 3. We end the paper with a couple of conjectures based on our results and some explicit computations  using Macaulay2 (Conjectures 5.1 and 5.3).

\vskip 4mm

\noindent \underline{Acknowledgment} Part of this work was done while the second author was a guest of the  Oslo and Akershus University College of Applied Sciences (which now is the Oslo Metropolitan University, Faculty of Technology, Art and Design)
 and she thanks the  Oslo Metropolitan  University  for its hospitality.

\vskip 2mm
\noindent \underline{Notation.} Throughout this paper $k$ will be an algebraically closed field of characteristic zero, $R=k[x_0, x_1, \cdots ,x_n]$, and $\mathfrak{m}=(x_0, \ldots,x_n)$.
For any graded Cohen-Macaulay quotient $A$ of $R$ of codimension $c$,
we let $I_A=\ker(R\twoheadrightarrow A)$ and  $K_A=\Ext^c_R (A,R)(-n-1)$ be its canonical module.  If $M$ is a
finitely generated graded $A$-module, let $\depth_{J}{M}$ denote
the length of a maximal $M$-sequence in a homogeneous ideal $J$
and let $\depth {M} = \depth_{\mathfrak m}{ M}$.


\section{Background and preparatory results}
For convenience of the reader we summarize the
 basic results on   determinantal ideals  and the associated complexes needed  in the sequel and we refer
to \cite{b-v}, \cite{eise}, \cite{L} and \cite{rm} for more details.

\begin{definition}
  \rm If $\cA$ is a homogeneous $q\times r$ matrix, we denote by $I(\cA)$ the
  ideal of $R$ generated by the maximal minors of $\cA$ and by $I_j(\cA)$ the
  ideal generated by the $j \times j$ minors of $\cA$. Assume $q\ge r$.
  $I(\cA)$ (resp. $R/I(\cA)$) is called a \emph{standard determinantal} ideal
  (resp. ring) if $\depth_{I(\cA)}R=q-r+1$, and $I_j(\cA)$ and
  $R/I_j(\cA)$ is said to be
  determinantal if $\depth_{I_j(\cA)} R=(q-j+1)(r-j+1)$.
 \end{definition}

Denote by
$\varphi
:F\longrightarrow G$ the morphism of free graded $R$-modules of
rank $t$ and $t+c-1$, defined by the homogeneous matrix $\cA$, and
 by ${\cC}_i(\varphi ^*)$ the (generalized) Koszul complex:
$${\cC}_i(\varphi ^*): \  0 \rightarrow \bigwedge^{i}G^* \otimes S_{0}(F^*)\rightarrow
\bigwedge^{i-1} G^*
\otimes S
_{1}( F^*)\rightarrow \ldots \rightarrow \bigwedge ^{0} G^* \otimes S_i (F^*) \rightarrow 0$$
where as usual $\bigwedge^{i}F^*$ (resp. $S_i(F^*)$) denotes the $i$-th exterior (symmetric) power of $F^*$.

 Let ${\cC}_i(\varphi^*)^*$ be the $R$-dual of
${\cC}_i(\varphi^*)$. The map $\varphi $ induces graded
morphisms $$\mu_i:\bigwedge ^{t+i}G^*\otimes
\bigwedge^tF\rightarrow \bigwedge^{i}G^*,$$
which are used to splice the complexes ${\cC}_{c-i-1}(\varphi^*)^*\otimes \bigwedge^{t+c-1}G^*\otimes \bigwedge^tF$ and
 ${\cC}_i(\varphi^*)$ to a complex ${\cD}_i(\varphi^*):$
 \begin{equation}\label{splice}
0 \rightarrow \bigwedge^{t+c-1}G^* \otimes S_{c-i-1}(F)\otimes
\bigwedge^tF\rightarrow \bigwedge^{t+c-2} G^* \otimes S _{c-i-2}(F)\otimes
\bigwedge ^tF\rightarrow \ldots \rightarrow \end{equation}
$$\bigwedge^{t+i}G^* \otimes
S_{0}(F)\otimes \bigwedge^tF\stackrel{\epsilon _{i}}{\longrightarrow} \bigwedge^{i} G^* \otimes S _{0}(F^*)
\rightarrow \bigwedge ^{i-1} G^* \otimes S_1(F^*)\rightarrow \ldots \rightarrow \bigwedge^0G^*\otimes
S_i(F^*)\rightarrow 0 .$$

The complex  ${\cD}_0(\varphi ^*)$ is called the {\em Eagon-Northcott
complex} and the complex  ${\cD}_1(\varphi ^*)$ is called the
{\em Buchsbaum-Rim complex}. Let us rename the complex  ${\cC}_c(\varphi ^*)$ as  ${\cD}_c(\varphi^*)$. Denote by $I_j(\varphi)$ the
ideal generated by the $j\times j$ minors of a matrix $\cA$
representing $\varphi$. Call $\varphi $ minimal if $\varphi (F)\subset \mathfrak{m}G$. Letting $S_{-1}M:=\Hom(M,R/I(\cA))$ we have the following well known result:

\begin{proposition}\label{resol} Let  $\varphi: F\rightarrow G$ be a graded minimal  morphism of free
$R$-modules of rank $t$ and $t+c-1$, respectively. Set $M= \coker(
\varphi ^*)$. Assume $\depth_{I_t(\varphi )} R=c$. Then it holds:

(i) ${\cD}_i(\varphi ^*)$ is acyclic for $-1\le i \le c$.

(ii) ${\cD}_0(\varphi ^*)$ is a minimal free graded $R$-resolution of
$R/I_t(\varphi )$  and  ${\cD}_i(\varphi ^*)$ is a minimal free graded
$R$-resolution of length $c$ of $S_iM$, $-1\le i \le c$ .

(iii) $K_{R/I_t(\varphi )} \cong S_{c-1}M(\mu )$ for some $\mu \in \ZZ$. So, up to twist,
$\mathcal{D}_{c-1}(\varphi ^*)$ is a minimal free graded $R$-module resolution of
$K_{R/I_t(\varphi )} $.

\end{proposition}

\begin{proof}
See, for instance \cite[Theorem 2.20]{b-v};   and \cite[Theorem A2.10 and Corollaries A2.12-13]{eise}.
\end{proof}

Let us also recall
 the following useful comparison of cohomology
groups. If $Z\subset X:=\Proj(B)$, $B=R/I_j(\varphi)$ is a closed subset such that $U=X\setminus Z$ is a local complete intersection,  $L$ and $N$ are finitely generated $B$-modules, $\widetilde{N}$ is locally free on
$U$ and
 $\depth_{I(Z)}L\ge r+1$, then the natural map
\begin{equation} \label{NM}
\Ext^{i}_{B}(N,L)\longrightarrow
\Hl_{*}^{i}(U,{\mathcal H}om_{\odi{X}}(\widetilde{N},\widetilde{L}))
\end{equation}
is an isomorphism, (resp. an injection) for $i<r$ (resp. $i=r$)
cf. \cite{SGA2}; expos\'{e} VI. Note that we interpret $I(Z)$ as $\mathfrak{m}$ if $Z=\emptyset $.


\section{Schur powers of the cokernel of a graded morphism}

Let $\varphi: F:=\oplus _{i=1}^tR(b_i)\longrightarrow G:=\oplus _{j=1}^{t+c-1} R(a_j)$ be a graded $R$-morphism between free $R$-modules of rank $t$ and $t+c-1$, respectively. Set $M=\coker (\varphi ^*)$.
It is a longstanding problem in commutative algebra to determine a graded free
minimal $R$-resolution of the $p$-exterior power, $\bigwedge^p M$, of $M$ or,
even more generally, a graded free minimal $R$-resolution of the cokernel of
the natural map
$$\varphi _{r,s}^{*}:\bigwedge ^rG^*\otimes \bigwedge ^sF^*\longrightarrow
\bigwedge ^{r+s}F^*, \quad f\otimes g \mapsto \bigwedge ^r\varphi ^*(f)\wedge
g.$$
Note that $\bigwedge ^pM=\coker (\varphi^*_{1,p-1})$ and, in particular,
$M=\coker(\varphi ^*_{1,0})$. So far, only few particular cases have been
described. To our knowledge, an explicit finite free $R$-resolution has been
computed for $\coker(\varphi ^*_{1,0})$ under the assumption that
$\depth_{I_t(\varphi)} R$ has the largest possible value, namely
$\depth_{I_t(\varphi)} R=c$. It is known that this value is
achieved if the entries of $\varphi $ are general enough. A minimal graded
free $R$-resolution of $\coker(\varphi ^*_{r,t-r})=R/I_r(\varphi )$ is known
under the assumption that
$\depth_{I_r(\varphi)} R=(t-r+1)(t-r+c)$ by \cite{L}. Again we
know that this value is achieved if the entries of $\varphi $ are general
enough. Finally, in \cite[Theorem 3.1 and Corollary 3.2]{BE75}, an explicit
finite free $R$-resolution of $\coker(\varphi ^*_{p,0})$ is given under the
hypothesis that $\depth_{I_t(\varphi )} R=c$, and in \cite[Proposition 6.1]{BE75} a minimal free $R$-resolution of $\bigwedge ^2M$ is
given in the case $c=1$ under the assumption
$\depth_{I_{t-1}(\varphi )} R=4$.

\vskip 2mm

In \cite[pg.\,299]{BE75}, Buchsbaum and Eisenbud predict the first terms of a graded free minimal $R$-resolution of $\coker (\varphi^*_{p-1,1})$ and, in particular, of $\wedge ^2 M$. To state their guess we need to fix some notation. Given a graded free $R$-module $F$, we define
\begin{equation}\label{LpqF}
L_p^qF=\ker (\delta_{p+1}^{q-1}) \ \ \text{ where } \\ \delta_k^{\ell }:S_{k-1} F\otimes \bigwedge ^{\ell }F\longrightarrow S_kF\otimes \bigwedge ^{\ell -1}F.
\end{equation}

\vskip 2mm
It holds (see, for instance, \cite[Corollary 2.3 and Proposition 2.5]{BE75}):
\begin{itemize}
\item[(1)] If $ p + q \ne  1$, then $L_p^q =\ker (\delta_{p+1}^{q-1})= \im (\delta_{p}^{q})=\coker (\delta_{p-1}^{q+1})$.
\item[(2)] $ L^{1}_pF = S_pF$ for all $p$.
\item[(3)] $L_1^qF = \bigwedge ^q F$ for all  $q \ne 0$.
\item[(4)] $L_p^0F = L_0^qF = 0$ for all $q \ne  1$ and all $p$.
\item[(5)] $L_p^qF = 0$ for all $q > \rk( F)$.
\item[(6)] If $\rk (F) = n$, then $L_p^nF=S_{p-1}F\otimes \bigwedge ^n F$.
\item[(7)]  If $\rk (F) = n$, then $\rk L_p^qF=\binom{n+p-1}{q+p-1}\binom{q+p-2}{p-1}$.
\end{itemize}

\vskip 2mm
Given the graded morphism  $\varphi ^*:G^*\longrightarrow F^*$  between free modules of rank $t+c-1$ and $t$, respectively, we set $$A^p(F):=\coker(\bigwedge ^{t-p} F\stackrel{\theta}{\longrightarrow} \bigwedge ^{t-p+1}F\otimes F^* )$$ where $\theta (\beta ):= \sum g _i\wedge \beta \otimes \gamma _i$ being $\{g_1,\cdots ,g_t\}$ is a basis of $F$ and $\{\gamma _1,\cdots ,\gamma_t\}$ its dual.

\begin{conj} \label{conjBE} (See \cite[pg. 299]{BE75}) \rm Let
  $\varphi ^*:G^*\longrightarrow F^*$ be a graded morphism between free
  $R$-modules of rank $t+c-1$ and $t$, respectively, and suppose
  $\depth_{I_{t-1}(\varphi)} R=2(c+1)$. A resolution of
  $\coker (\varphi^*_{p-1,1})$ should always start out as:
$$ \cdots \longrightarrow L_2^{p-1}(G)^* \oplus A^p(F)\otimes \bigwedge ^tG^*  \longrightarrow \bigwedge ^{p-1}G^* \otimes F^* \longrightarrow \bigwedge ^p F^*\longrightarrow  \coker (\varphi^*_{p-1,1}) \longrightarrow 0.$$
In particular, a resolution of $\bigwedge ^2M$ should always start out as:
$$ \cdots \longrightarrow L_2^{1}(G)^* \oplus A^2(F)\otimes \bigwedge ^tG^*  \longrightarrow G^* \otimes F^* \longrightarrow \bigwedge ^2 F^*\longrightarrow \bigwedge ^2M \longrightarrow 0.$$
(We will see later that $L_2^{1}(G)^* \oplus A^2(F)\otimes \bigwedge ^tG^*\cong S_2(G)^* \oplus S_2(F^*)\otimes \bigwedge ^t F\otimes \bigwedge ^tG^*$.)\end{conj}

\begin{remark} \rm The assumption $\depth_{I_{t-1}(\varphi)} R=2(c+1)$ cannot
  be skipped. Indeed we have used Macaulay2 to compute the minimal resolution
  of $\wedge ^2M$ for general linear $3 \times (c+2)$ matrices with
  coefficients in $R={\mathbb Q}[x_0, x_1, \cdots ,x_n]$ for different values
  of $n $ and any $c \in \{2,3\}$ with $\mathbb Q$ the field of rational
  numbers. The graded Betti numbers for fixed $c$ and various $n \ge 2(c+1)-1$
  are the same and satisfies Conjecture~\ref{conjBE} while for the choice
  $n+1 = 2(c+1)-1$ the minimal resolution of $\wedge ^2M$ radically change and
  the prediction in Conjecture~\ref{conjBE} does not hold. Note that in the
  latter case, $\depth_{I_{t-1}(\varphi)} R=2(c+1)-1$ and
  $R/ I_{t-1}(\varphi)$ is not determinantal.
\end{remark}
\vskip 2mm

The main goal of this section is to relate the canonical module of
$I_{i}(\varphi)$ with a suitable Schur power of $M=\coker (\varphi ^*)$
(Proposition \ref{canonical module} and Theorem \ref{main2}) and prove
Conjecture \ref{conjBE} for $p=2$ and $c=2$ (Theorems \ref{BE's_guess}).
Indeed, we will determine a minimal graded free resolution of
$\bigwedge ^p M=\coker (\varphi^*_{1,p-1})$ for $p\ge 1$ and $c=2$ (Theorem
\ref{BE's_guess}). Our proof is based on the fact that the minimal graded free
$R$-resolution of the determinantal ideal $I_{t+1-p}(\varphi )$ is known, for
any $c\ge 2$, under the assumption that
$\depth_{I_{t-p+1}(\varphi)} R=p(p+c-1)$, and Proposition
\ref{canonical module} which identifies, up to twist, $\bigwedge ^pM$ with the
canonical module of $R/I_{t+1-p}(\varphi ) $ provided $c=2$. In section 4, we
will generalize Theorem \ref{BE's_guess} to the case $c=3$ and we will prove
Conjecture \ref{conjBE} for $p=2$ and $c=3$ up to an explicitly described
summand which we guess splits off (Theorem \ref{resolwedge2}).

For sake of completeness, let us start  explaining how we build a minimal  graded free $R$-resolution of the  determinantal ideal $I_{t+1-p}(\varphi )$.
 This construction rests heavily on the theory of Schur functors and we refer to \cite{L} and \cite{Bu} for more details. We know that $R/I_i(\varphi )$ is a Cohen-Macaulay ring of projective dimension $(t-i+1)(t+c-i)$ and let  us call $\LL^{\bullet}(i,\varphi)$ a graded minimal free $R$-resolution of $R/I_i(\varphi )$.
 To define the  terms $$\LL^{-s}(i,\varphi) \ \text{ for } \  1\le s \le  (t-i+1)(t+c-i)$$
  of the complex $\LL^{\bullet}(i,\varphi)$,
 we need to recall some basic facts on partitions and Schur functors.

 \vskip 2mm
\noindent  {\bf Partitions.}
 Let $r$ be an integer. A {\em partition} is an $r$-tuple $I=(i_1,\cdots , i_r)\in \ZZ^r_{\ge 0}$ such that $0\le i_1 \leq i_2 \leq
 \cdots \leq i_r$.  The {\em length} of the partition $I$ is the number of its non-zero elements, and its {\em weight}, denoted by $|I|$, is $i_1+\cdots +i_r$.
Given two partitions $I=(i_1,\cdots , i_r)\in \ZZ_{\ge 0}^r$ and $I=(j_1,\cdots , j_s)\in \ZZ_{\ge 0}^s$ we say that $J\subset I$ if $j_v\le i_v$ for all $v\ge 1$.

\noindent We will draw the {\em Young diagram} (also called  {\em Ferrer diagram}) corresponding to a partition $I$ by
putting $i_j$ boxes in the $j$-th row and to the left as in the pictures below. We will also denote  by $I$ the
corresponding Young diagram. If we flip this diagram over its  diagonal, we get the Young diagram of the {\em conjugate} (or {\em dual}) partition $I^*$. For example the partition $I=(2,2,3,5,8)$ corresponds
to the Young diagram
\vspace{5mm}
\smallbreak
\centerline{
\vbox{
\def\star{\rlap{\hbox to 13pt{\hfil\raise3.5pt\hbox{$*$}\hfil}}}
\def\ {\hbox to 13pt{\vbox to 13pt{}\hfil}}
\def\*{\star\ }
\def\_{\hbox to 13pt{\hskip-.2pt\vrule\hss\vbox to 13pt{\vskip-.2pt
            \hrule width 13.4pt\vfill\hrule\vskip-.2pt}\hss\vrule\hskip-.2pt}}
\def\x{\star\_}
\offinterlineskip
\hbox{\ \ \     \_\_\ \ \ \ \ \ \ \ \ }
\hbox{\ \ \     \_\_\ \ \ \ \ \ \ \ \ }
\hbox{\ \ \    \_\_\_\ \ \ \ \ \ \ \ }
\hbox{\ \ \    \_\_\_\_\_\ \ \ \ \ \ }
\hbox{\ \ \    \_\_\_\_\_\_\_\_\ \ \ }
\hbox{\ \ \   \  \  \  \  \  \  \  \  \ \ \ }
}
}
\noindent and $I^*=(1,1,1,2,2,3,5,5)$ corresponds to the Young diagram:
 \vspace{5mm}
\smallbreak
\centerline{
\vbox{
\def\star{\rlap{\hbox to 13pt{\hfil\raise3.5pt\hbox{$*$}\hfil}}}
\def\ {\hbox to 13pt{\vbox to 13pt{}\hfil}}
\def\*{\star\ }
\def\_{\hbox to 13pt{\hskip-.2pt\vrule\hss\vbox to 13pt{\vskip-.2pt
            \hrule width 13.4pt\vfill\hrule\vskip-.2pt}\hss\vrule\hskip-.2pt}}
\def\x{\star\_}
\offinterlineskip
\hbox{\ \ \     \_\ \ \ \ \ \ \ \ \ \ }
\hbox{\ \ \     \_\ \ \ \ \ \ \ \ \ \ }
\hbox{\ \ \     \_\ \ \ \ \ \ \ \ \ \ }
\hbox{\ \ \     \_\_\ \ \ \ \ \ \ \ \ }
\hbox{\ \ \     \_\_\ \ \ \ \ \ \ \ \ }
\hbox{\ \ \    \_\_\_\ \ \ \ \ \ \ \ }
\hbox{\ \ \    \_\_\_\_\_\ \ \ \ \ \ }
\hbox{\ \ \    \_\_\_\_\_\ \ \ \ \ \ }
}
}

\vskip 4mm
\noindent  Given a partition $I=(i_1,\cdots , i_r)$, we will say that $I$ has {\em Durfee square} $p_I$ if $p_I$ is the dimension of the greatest square contained in the Young diagram of $I$. So, if $I=(2,2,3,5,8)$, then $p_I=3$,

 \vskip 3mm
\noindent  {\bf Schur functors.}
Given a rank $r$ vector bundle $\shE$ on $\PP^n$, we will denote by
$\Sigma^{I}\shE$  the result of applying the Schur functor $\Sigma^{I}$ with
 $I=(i_1, \cdots, i_s) \in \ZZ^{s}_{\ge 0}$  a partition with $s\le r$ parts, i.e., $0\le i_1 \leq i_2 \leq
 \cdots \leq i_s$. As special cases one recover symmetric and exterior powers:
\[ \Sigma^{(0,\cdots,0,m)}\shE = S_{m}\shE  \]
is the $m$-symmetric power of $\shE$ and
\[ \Sigma^{(0, \cdots,0,\tiny{\overbrace{1,\cdots,1}^{m}})}\shE = \bigwedge^m\shE  \]
is the $m$-exterior power of $\shE$.  We also use the convention $\Sigma^{(i_1, \cdots,
i_s)}\shE=  \Sigma^{(\tiny{\overbrace{0,\cdots ,0,i_1, \cdots, i_s}^{r}})} \shE $ whenever $s<r$.
 In particular, for any rank $r$ vector bundle $\shE$ on $\PP^n$,
\[ \det(\shE)= \bigwedge^r\shE=   \Sigma^{(\tiny{\overbrace{1,\cdots,1}^{r}})}\shE, \text{ and }  \]
\[ \Sigma^{(i_1+m, \cdots,
i_r+m)}\shE=\Sigma^{(i_1, \cdots, i_r)}\shE \otimes  \det(\shE)^{\otimes m}. \]

In general, given  $I=(i_1, \cdots, i_q)$  a partition and $I^*=(j_1, \cdots, j_s)$ its dual
and writing $(\otimes ^{j_1}\shE )\otimes \cdots \otimes (\otimes ^{j_s} \shE)$ as $ \bigotimes ^{|I|}(\shE)$, we define the {\em Schur power} $\Sigma ^I\shE$  of $\shE$ with exponent $I$ as the image of the composite map
$$d_I(\shE):\bigwedge^{j_1}\shE \otimes \cdots \otimes \bigwedge^{j_s} \shE\stackrel{\alpha _I}{\longrightarrow} \bigotimes ^{|I|}(\shE)\stackrel{\beta _I}{\longrightarrow} S_{i_1}\shE\otimes \cdots \otimes S_{i_q}\shE$$
where $\alpha _I$ is the tensor product of the injections $\bigwedge^{j_v} \shE \longrightarrow \otimes ^{j_v}\shE$ which sends
$f_v^1\wedge \cdots \wedge f_v^{j_v}\mapsto \frac{1}{j_v!}\sum _{\sigma \in \mathfrak{S}_{j_v} }(-1)^{\varepsilon(\sigma)}f_v^{\sigma(1)}\otimes \cdots \otimes f_v^{\sigma(j_v)}$ being $\mathfrak{S}_{j_v}$ the permutation group. Moreover,
$\beta _I:\bigotimes ^{|I|}(\shE) \longrightarrow S_{i_1}\shE\otimes \cdots
\otimes S_{i_q}\shE$ sends $(f_1^1\otimes \cdots \otimes f_1^{j_1})\otimes
(f_2^1\otimes \cdots \otimes f_2^{j_2})\otimes \cdots \otimes (f_s^1\otimes
\cdots \otimes f_s^{j_s})\mapsto (f_1^qf_2^q \cdots f_s^q)\otimes
(f_1^{q-1}f_2^{q-1} \cdots f_s^{q-1})\otimes \cdots \otimes (f_1^1f_2^1 \cdots f_s^1) $ where  $q=j_s$ and where we let $f_{i}^{j}=1$ for $j_i<j\le j_s$, i.e. $f_j^i=1$ if $(i,j)$ is not a position
(or box)
 in the Young diagram of $I^*$.
Recall that if $\shE$ is a rank $r$ vector bundle on $\PP^n$ and $I=(i_1,\cdots ,i_r)$ a partition ($0\le i_1\le i_2\le \cdots \le i_r$) then (see, for example, \cite[Theorem 6.3]{FH}):
\begin{equation} \label{rank}
\rk (\Sigma ^I \shE )= \prod _{1\le u<v\le r} \frac{i_v-i_u+v-u}{v-u}.
\end{equation}

For instance,  we consider $I=(1,2,4)$ and its dual $I^*=(1,1,2,3)$ with Young diagrams

\vspace{5mm}
\smallbreak
\centerline{
\vbox{
\def\star{\rlap{\hbox to 13pt{\hfil\raise3.5pt\hbox{$*$}\hfil}}}
\def\ {\hbox to 13pt{\vbox to 13pt{}\hfil}}
\def\*{\star\ }
\def\_{\hbox to 13pt{\hskip-.2pt\vrule\hss\vbox to 13pt{\vskip-.2pt
            \hrule width 13.4pt\vfill\hrule\vskip-.2pt}\hss\vrule\hskip-.2pt}}
\def\x{\star\_}
\offinterlineskip
\hbox{\ \ \     \_\ \ \ \ \ \ \ \ \ \ }
\hbox{\ \ \    \_\_\ \ \ \ \ \ \ \ \ }
\hbox{\ \ \    \_\_\_\_\ \ \ \ \ \ \ }
}
\vbox{
\def\star{\rlap{\hbox to 13pt{\hfil\raise3.5pt\hbox{$*$}\hfil}}}
\def\ {\hbox to 13pt{\vbox to 13pt{}\hfil}}
\def\*{\star\ }
\def\_{\hbox to 13pt{\hskip-.2pt\vrule\hss\vbox to 13pt{\vskip-.2pt
            \hrule width 13.4pt\vfill\hrule\vskip-.2pt}\hss\vrule\hskip-.2pt}}
\def\x{\star\_}
\offinterlineskip
\hbox{\ \ \     \_\ \ \ \ \ \ \ \ \ \ }
\hbox{\ \ \     \_\ \ \ \ \ \ \ \ \ \ }
\hbox{\ \ \    \_\_\ \ \ \ \ \ \ \ \ }
\hbox{\ \ \    \_\_\_\ \ \ \ \ \ \ \ }
}
}
\vskip 4mm \noindent Then,  $\Sigma ^I\shE$  is the image of the  map

$$\hskip -5.5cm d_I(\shE):\shE \otimes \shE \otimes \bigwedge ^2\shE \otimes \bigwedge^{3}\shE \longrightarrow \shE \otimes  S_{2}\shE \otimes  S_{4}\shE$$
 $$  f_1^1\otimes f_2^1 \otimes (f_3^1\wedge f_3^2)\otimes (f_4^1\wedge f_4^2\wedge f_4^3)\mapsto \frac{1}{12}
 \sum _{\overset{\tau \in \mathfrak{S}_2}{\sigma \in \mathfrak{S}_3}} (-1)^{\varepsilon(\tau)+ \varepsilon(\sigma)} f_4^{\sigma(3)}\otimes (f_3^{\tau(2)}f_4^{\sigma(2)})\otimes (f_1^1f_2^1f_3^{\tau(1)}f_4^{\sigma(1)})$$
 obtained composing the maps

 $$\hskip -7.8cm \alpha _I : \shE \otimes \shE \otimes \bigwedge ^2\shE \otimes \bigwedge^{3}\shE \longrightarrow  \bigotimes ^{7}\shE$$ $$ f_1^1\otimes f_2^1 \otimes (f_3^1\wedge f_3^2)\otimes (f_4^1\wedge f_4^2\wedge f_4^3)\mapsto \frac{1}{12}\sum _{\overset{\tau \in \mathfrak{S}_2}{ \sigma \in \mathfrak{S}_3}} (-1)^{\varepsilon(\tau)+\varepsilon(\sigma)}f_1^1\otimes f_2^1 \otimes
 (f_3^{\tau(1)}\otimes f_3^{\tau(2)}) \otimes
  (f_4^{\sigma(1)}\otimes f_4^{\sigma(2)}\otimes f_4^{\sigma(3)})$$
and

$$\hskip 2cm \beta _I:\bigotimes ^{7}(\shE) \longrightarrow \shE \otimes  S_{2}\shE \otimes  S_{4}\shE$$ $$  f_1^1\otimes f_2^1 \otimes (f_3^1\otimes f_3^2)\otimes (f_4^1\otimes f_4^2\otimes f_4^3) \mapsto
 f_4^{3}\otimes (f_3^{2}f_4^{2})\otimes (f_1^1f_2^1f_3^{1}f_4^{1}).$$
In addition, if $\shE$ has rank 5, then $\Sigma ^I \shE$ has rank 700.

\vskip 4mm
We will often use {\em Pieri's formula} (see \cite[pg. 209]{L}): Let $E$ be a free $R$-module of rank $r$ and let $I$ be a partition. It holds:
\begin{equation}
\label{Pieri1} \wedge ^nE\otimes \Sigma ^IE\cong \bigoplus _J\Sigma ^J E \end{equation}
where the index set consists of all partitions $J$ whose Young diagram is
obtained by adding n boxes to the Young diagram of $I$ with no two boxes in
the same row, see \cite[Corollary IV.2.6]{ABW}.
So, if $I=(i_1,\cdots, i_s)$ then the index set consists of all
$J=(j_0,j_1,\cdots,j_{s})$ satisfying $|J|=|I|+n$, $i_t\le j_t\le i_t+1$ and
$j_0=0$ (i.e. $J=(j_1,\cdots,j_{s})$) if $s=r$ and $0 \le j_0\le 1$
if $s < r$.

Dually, we have
\begin{equation}
\label{Pieri2}  S_nE\otimes \Sigma ^IE\cong \bigoplus _J\Sigma ^J E \end{equation}
where the index set consists of all partitions $J$ whose Young diagram is obtained by adding n boxes to the Young diagram of $I$ with no two boxes in the same column.

\begin{example} \rm
Let $E$ be a free $R$-module of rank 6. By Pieri's formula we have
$$S_4E\otimes \bigwedge ^2E\cong \Sigma ^{(0,0,0,0,1,5)}E\oplus \Sigma ^{(0,0,0,1,1,4)}E,$$
i.e. $S_4E\otimes \bigwedge ^2E$ is a rank 1890 free $R$-module which decompose into a direct sum of two free $R$-modules of rank 1050 and 840, respectively.
\end{example}

\vskip 2mm

Let  $\varphi: F\longrightarrow G$ be a graded $R$-morphism of free $R$-modules of rank $t$ and $t+c-1$, respectively.  We define the complex $\LL^{\bullet}(i,\varphi)$ of $R$-modules as follows:

\begin{equation}\label{Lj} L^{j}(i,\varphi)= \bigoplus _{\overset{-|I|+n(I)=j}{ I,I'\ne \emptyset }}\Sigma ^{(I)^*}(G^*)\otimes \Sigma ^{I'}(F), \quad \text{ for } \quad  -(t-i+1)(t+c-i)\le j \le -1,\end{equation}

\noindent where  $I=(i_1, \cdots, i_q)$  is a partition, $q=t-i+1$, $p:=p_I$ is the Durfee  square of $I$, $n(I)=(i-1)p_I$, $I'=(i_1, \cdots ,i_{q-p_{I}},\overbrace{p_{I},\cdots, p_{I}}^{i-1},i_{q-p_{I}+1}-i+1,\cdots, i_q-i+1)$ if $i_{q-p_{I}+1}\ge p_I+i-1$, otherwise $I'=\emptyset $ and $(I)^*$ the dual of the partition $I$.

Now, we have to define a morphism between  $ L^{j}(i,\varphi)$ and $ L^{j+1}(i,\varphi)$, i.e. for any pair of partitions $(I,H)$ with $-|I|+n(I)=j=-|H|+n(H)-1$ we must define $$\psi _{I,H}: \Sigma ^{(I)^*}(G^*)\otimes \Sigma ^{ I'}(F)  \longrightarrow \Sigma ^{(H)^* }(G^*)\otimes \Sigma ^{ H'}(F).$$ We distinguish two cases:

Case 1: If   $H\nsubseteq I$, we set $\psi _{(I,H)}=0$.

Case 2: $H\subset I$. By \cite[Lemma 5.13]{L}, if $(I,H)$ is a pair of partitions such that $H\subset I$, $I'\ne \emptyset $, $H'\ne \emptyset$ and $-|I|+n(I)=j=-|H|+n(H)-1$, then $I$ and $H$ (resp. $I'$ and $H'$) coincide except in one spot, i.e.
$I=(i_1, \cdots , i_{x-1},a,i_{x+1},\cdots ,i_q)$, $H=(i_1, \cdots , i_{x-1},h,i_{x+1},\cdots ,i_q)$, $(I')^*=(j_1, \cdots , j_{y-1},a',j_{y+1},\cdots ,j_s)$, $(H')^*=(j_1, \cdots , j_{y-1},h',j_{y+1},\cdots ,j_s)$. By \cite[Lemma 5.13]{L}  (see also \cite[Lemma 3.1]{Bu}), we have $a-h=a'-h'$ and we set  $\rho : =a-h=a'-h'$.
 Since the Young diagram of $I^*$ is
obtained by adding $\rho$ boxes to the Young diagram of $H^*$ with no two
boxes in the same row, we have by Pieri's formula (\ref{Pieri1}) that $\Sigma ^{(I)^*}G^*$ is a direct summand of $ \bigwedge ^{\rho}G^*\otimes \Sigma ^{(H)^*} G^*$. Therefore, we have a canonical injective map
  $$\delta _{I,\rho}:\Sigma ^{(I)^*}G^*\longrightarrow \bigwedge ^{\rho}G^*\otimes \Sigma ^{(H)^*} G^* $$
    and, dually, we have a canonical injective map
     $$\delta _{I',\rho}:\Sigma ^{I'}F\longrightarrow  \bigwedge ^{\rho}F\otimes \Sigma ^{H'} F$$
which give us
 $$\delta _{I,\rho}\otimes \delta _{I',\rho}:\Sigma ^{(I)^*}G^*\otimes \Sigma ^{I'}F\longrightarrow \bigwedge ^{\rho}G^*\otimes \bigwedge ^{\rho}F \otimes \Sigma ^{(H)^*} G^*\otimes \Sigma ^{H'} F.$$

\noindent Composing the  morphism $\delta _{I,\rho}\otimes \delta _{I',\rho}$
with $$\mu \otimes 1:\bigwedge ^{\rho}G^*\otimes \bigwedge ^{\rho}F \otimes
\Sigma ^{(H)^*} G^*\otimes \Sigma ^{H'} F\longrightarrow  \Sigma ^{(H)^*}
G^*\otimes \Sigma ^{H'} F$$ where $\mu: \bigwedge ^{\rho}G^*\otimes \bigwedge
^{\rho}F \longrightarrow R$ is the natural map induced by $\varphi ^*:G^*\longrightarrow F^*$, we get
 a morphism
 $$\sigma _{(I,H)}:\Sigma ^{(I)^*}G^*\otimes \Sigma ^{
 I'}F\longrightarrow  \Sigma ^{(H)^* } G^*\otimes \Sigma ^{H'} F.$$
%
%
%
 Finally, we define the differential map:
 $$d_j=\oplus \psi _{(i,H)}: L^{j}(i,\varphi)=\bigoplus _{\overset{-|I|+n(I)=j}{ I,I'\ne \emptyset}} \Sigma ^{(I)^*}(G^*)\otimes \Sigma ^{ I'}(F)\longrightarrow  L^{j+1}(i,\varphi)=
 \bigoplus _{\overset{-|H|+n(H)=j+1}{ H,H'\ne \emptyset }}\Sigma ^{(H)^*}(G^*)\otimes \Sigma ^{ H'}(F).$$

We have:

\begin{theorem}\label{resoldet} Let $\varphi: F\longrightarrow G$ be a graded
  $R$-morphism of free $R$-modules of rank $t$ and $t+c-1$, respectively.
  Assume that $\varphi $ is minimal and
  $\depth_{I_{i}(\varphi)} R=(t-i+1)(t+c-i)$. Then, the
  complex $\LL^{\bullet}(i,\varphi)$ is exact and
  $\Hl ^0(\LL^{\bullet}(i,\varphi))=R/I_{i}(\varphi )$, i.e.
  $\LL^{\bullet}(i,\varphi)$ is a graded minimal free $R$-resolution of
  $R/I_{i}(\varphi ).$
\end{theorem}

\begin{proof} See \cite[Th\'{e}or\`{e}me 3.3]{L} or \cite[Theorem 3.2]{Bu}.
\end{proof}

Our next goal will be to apply the above Theorem to solve Conjecture \ref{conjBE} for $p=c=2$ (see Theorem \ref{BE's_guess}).

\begin{proposition}\label{canonical module} Set $R=k[x_0,\cdots ,x_n]$. Let $\varphi ^*:G^*\longrightarrow F^*$ be a graded minimal morphism between free $R$-modules of rank $t+1$ and $t$, respectively.  Set $B_{t-p+1}:=R/I_{t+1-p}(\varphi )$, $M=\coker (\varphi ^*)$ and $\ell :=\sum _{j=1}^{t+1}a_j-\sum _{i=1}^tb_i$ and assume that $B_{t+1-p}$ is determinantal. Then, we have:
$$K_{B_{t+1-p}}(n+1-p\ell) \cong  \bigwedge ^p M .$$
\end{proposition}

\begin{proof} By hypothesis $M=\coker(\varphi ^*)$. Therefore, a minimal graded $R$-free resolution of $\bigwedge ^pM$  starts out as:
$$ \cdots  \longrightarrow G^* \otimes \bigwedge ^{p-1}F^*  \stackrel{\overline{\varphi ^*}}{\longrightarrow} \bigwedge ^p F^*\longrightarrow \bigwedge ^pM \longrightarrow 0$$ with $\overline{\varphi ^*}(\eta \otimes \rho)=\varphi ^*(\eta)\wedge \rho$. So, $\bigwedge ^pM=\coker(\overline{\varphi ^*})$.

On the other hand, we know that $B_{t-p+1}$ is a Cohen-Macaulay ring and $\LL^{\bullet}(t-p+1,\varphi)$ is a graded minimal free $R$-resolution of $B_{t-p+1}$ of length $p(p+1)$ (see Theorem \ref{resoldet}). Applying $\Hom_R(-,R)$ to  $\LL^{\bullet}(t+1-p,\varphi)$ we get  a minimal graded free $R$-resolution of $\Ext_R^{p(p+1)}(B_{t-p+1},R)\cong K_{B_{t+1-p}}(n+1)$. Let us explicitly compute  $\LL^{-p(p+1)}(t+1-p,\varphi)$, $\LL^{-p(p+1)+1}(t+1-p,\varphi)$ and the morphism linking them.
 $\LL^{-p(p+1)}(t+1-p,\varphi)$ corresponds to the unique partition $I=(\overbrace{t+1, \cdots , t+1}^p)$ for which the Durfee square $p_I=p$, $n(I)=p(t-p)$ and $I'=(\overbrace{p,\cdots ,p}^{t-p},\overbrace{p+1,\cdots ,p+1}^p)$. Thus, we have $$\begin{array}{rcl} \LL^{-p(p+1)}(t+1-p,\varphi) & := & \Sigma ^{(I)^*}G^*\otimes \Sigma ^{I'}F \\
 & \cong & (\bigwedge ^{t+1} G^*)^{\otimes p}\otimes (\bigwedge ^tF)^{\otimes p}\otimes \bigwedge ^p F \\
 & \cong & \bigwedge ^p F\otimes R(-p\ell ).\end{array}$$

\noindent Analogously, we apply \cite[Corollaire 5.10]{L} and we get that the $R$-free module
$\LL^{-p(p+1)+1}(t+1-p,\varphi)$ corresponds to a unique partition, namely,  $I=(t,\overbrace{t+1, \cdots , t+1}^{p-1})$ for which $p_I=p$, $n(I)=p(t-p)$, $I'=(\overbrace{p,\cdots ,p}^{t-p+1},\overbrace{p+1,\cdots ,p+1}^{p-1})$ and $I^*=(p-1, \overbrace{p,\cdots ,p}^{t})$. Thus, we have $$\begin{array}{rcl} \LL^{-p(p+1)+1}(t+1-p,\varphi)& := & \Sigma ^{(I)^*}G^*\otimes \Sigma ^{I'}F \\ & \cong & (\bigwedge ^{t+1} G^*)^{\otimes (p-1)}\otimes \bigwedge ^t(G^*)\otimes (\bigwedge ^tF)^{\otimes p}\otimes \bigwedge ^{p-1} F \\
& \cong & G\otimes \bigwedge ^{p-1} F\otimes R(-p\ell ) \end{array}$$
where the last isomorphism follows from the well known fact that $\bigwedge ^t(G^*)\cong G\otimes \bigwedge ^{t+1}(G^*)$.
To  finish the proof it is enough to observe that for the unique partitions $I=(\overbrace{t+1, \cdots , t+1}^p)$ and $H=(t,\overbrace{t+1, \cdots , t}^{p-1})$ associated to $\LL^{-p(p+1)}(t+1-p,\varphi)$ and $\LL^{-p(p+1)+1}(t+1-p,\varphi)$, respectively,
we have $\psi _{(I,H)}^* \cong \overline{\varphi ^*}\otimes 1_{R(p\ell)}$, whence  $K_{B_{t+1-p}}(n+1-p\ell) \cong  \bigwedge ^p(M) .$
\end{proof}

By \cite[Proposition 4]{V}, we know that $\bigwedge ^pM$ is a maximal Cohen-Macaulay $R/I_{t+1-p}(\varphi)$-module. As a consequence of the above Proposition we recover Vetter's result and even more. We have:

\begin{corollary} \label{gorensteinessExteriorpowerM} Set
  $R=k[x_0,\cdots ,x_n]$. Let $\varphi ^*:G^*\longrightarrow F^*$ be a graded
  minimal morphism between free $R$-modules of rank $t+1$ and $t$,
  respectively, and let $M$ be the cokernel of $\varphi ^*$. Assume that
  $\depth_{I_{t+1-p}(\varphi)} R=(p+1)p$. Then,
  $\bigwedge ^p M $ is a Gorenstein $R/I_{t+1-p}(\varphi)$-module.
\end{corollary}

\begin{proof} Since $R/I_{t+1-p}(\varphi)$ is a Cohen-Macaulay ring, the result immediately follows from Proposition \ref{canonical module} and \cite[Lemma 5.4]{KMMNP}.
\end{proof}

\begin{theorem} \label{BE's_guess} Set $R=k[x_0,\cdots ,x_n]$. Let
  $\varphi ^*:G^*\longrightarrow F^*$ be a graded minimal morphism between
  free $R$-modules of rank $t+1$ and $t$, respectively. Set
  $\ell :=\sum _{j=1}^{t+1}a_j-\sum _{i=1}^tb_i$ and assume
  $\depth_{I_{t+1-p}(\varphi)} R=p(p+1)$. Then, the $R$-dual
  $\LL^{\bullet}(t+1-p,\varphi)^{*}$ of $\LL^{\bullet}(t+1-p,\varphi)$ is a
  minimal free $R$-resolution of $\bigwedge ^p M (p\ell )$.

  In particular, Conjecture \ref{conjBE} (with $ \varphi $ minimal) is true
  for $p=c=2$.
\end{theorem}

\begin{proof} We know that $\LL^{\bullet}(t+1-p,\varphi)$ is a minimal graded free $R$-resolution of $B_{t+1-p}:=R/I_{t+1-p}(\varphi)$. We also know that $R/I_{t+1-p}(\varphi)$ is a Cohen-Macaulay ring. Therefore, the $R$-dual $\LL^{\bullet}(t+1-p,\varphi)^{*}$ of $\LL^{\bullet}(t+1-p,\varphi)$ gives us a minimal graded free $R$-resolution of $K_{B_{t+1-p}(\varphi)}(n+1)$ and we now conclude using Proposition \ref{canonical module}.

To prove Conjecture \ref{conjBE}  for $p=c=2$, we only have to identify the last terms of  $\LL^{\bullet}(t+1-p,\varphi)$. By the proof of Proposition \ref{canonical module} we know that $$\LL^{-p(p+1)}(t+1-p,\varphi)=(\bigwedge ^{t+1} G^*)^p\otimes (\bigwedge ^tF)^p\otimes \bigwedge ^p F\cong \bigwedge ^p F(-p\ell )$$ and $$\LL^{-p(p+1)+1}(t+1-p,\varphi)=(\bigwedge ^{t+1} G^*)^p\otimes G \otimes (\bigwedge ^tF)^p\otimes \bigwedge ^{p-1} F\cong G\otimes \bigwedge ^{p-1} F(-p\ell ).$$ Let us compute  $\LL^{-p(p+1)+2}(t+1-p,\varphi)$.
$\LL^{-p(p+1)+2}(t+1-p,\varphi)$ corresponds to two different partitions:
\begin{itemize}
\item $I_1:=(p-1,\overbrace{t+1, \cdots , t+1}^{p-1})$ with Durfee square $p_{I_1}=p-1$,  $n(I_1)=(p-1)(t-p)$ and $I_1'=( \overbrace{p-1,\cdots ,p-1}^{t-p+1},\overbrace{p+1,\cdots ,p+1}^{p-1})$.
 \item $I_2:=(t,t,\overbrace{t+1, \cdots , t+1}^{p-2})$ for which $p_{I_2}=p$, $n(I_2)=p(t-p)$ and $I_2'=( \overbrace{p,\cdots ,p}^{t-p+2},\overbrace{p+1,\cdots ,p+1}^{p-2})$
\end{itemize}
Therefore, we have
 $$\begin{array}{rcl} \Sigma ^{I_1'}F & = &  (\bigwedge ^t F)^{\otimes p-1}\otimes \Sigma ^{(0,\cdots,0,\overbrace{2,\cdots ,2}^{p-1})}F, \\
 \Sigma ^{I_2'}F &  = &  (\bigwedge ^tF)^{\otimes p}\otimes \Sigma ^{(0,\cdots,0,\overbrace{1,\cdots ,1}^{p-2})}F \\
 & \cong & (\bigwedge ^tF)^{\otimes p}\otimes \bigwedge ^{p-2}F, \end{array}
 $$
 $$\begin{array}{rcl}
  \Sigma ^{(I_1)^* }G^* & = &  (\bigwedge ^{t+1} G^*)^{\otimes p-1} \otimes \Sigma ^{(0,\cdots,0,\overbrace{1,\cdots ,1}^{p-1})}G^* \\
  &  \cong &   (\bigwedge ^{t+1} G^*)^{\otimes p-1}\otimes \bigwedge ^{p-1}G^* , \text{ and } \\
 \Sigma ^{(I_2)^* }G^* & = & (\bigwedge ^{t+1}G^*)^{\otimes p-2}\otimes  \Sigma ^{(0,2,\cdots ,2)}G^* \\
 & \cong & (\bigwedge ^{t+1}G^*)^{\otimes p}\otimes  \Sigma ^{(-2,0,\cdots ,0)}G^* \\
 & \cong & (\bigwedge ^{t+1}G^*)^{\otimes p}\otimes S_2G.  \end{array}
 $$
In particular, for $p=c=2$, we obtain:
 $$\begin{array}{rcl} \Sigma ^{I_1'}F & = &  \bigwedge ^t F\otimes S_2F  \cong   (\bigwedge ^tF)^{\otimes 2}\otimes \bigwedge ^t F^*\otimes S_2F, \\
 \Sigma ^{I_2'}F &  = &  (\bigwedge ^tF)^{\otimes 2}, \\
  \Sigma ^{(I_1)^* }G^* & = &  \bigwedge ^{t+1} G^* \otimes G^*  \cong   (\bigwedge ^{t+1} G^*)^{\otimes 2}\otimes \bigwedge ^tG , \text{ and } \\
 \Sigma ^{(I_2)^* }G^* & = & (\bigwedge ^{t+1}G^*)^{\otimes 2}\otimes S_2G.  \end{array}
 $$
 Thus, we have $$\begin{array}{rcl} \LL^{-p(p+1)+2}(t+1-p,\varphi) & := & ( \Sigma ^{(I_1)^* }G^*\otimes \Sigma ^{I_1'}F)\oplus (\Sigma ^{(I_2)}G^*\otimes \Sigma ^{I_2'}F )\\
  & \cong & \big( S_2G\oplus \bigwedge ^tG\otimes S_2F\otimes \bigwedge ^{t} F^*\big)\otimes R(-2\ell ).\end{array}$$
Since, for $c=p=2$, we have $A^2(F)=\coker(\bigwedge ^{t-2}F\longrightarrow \bigwedge^{t-1}F\otimes F^*)\cong S_2F^*\otimes \bigwedge ^tF$, we have got what we wanted to prove.
\end{proof}

\begin{remark}\rm (1) The above Theorem solves Conjecture \ref{conjBE}  for
  $p=c=2$ and it gives a minimal resolution of $\coker (\varphi ^* _{1,p-1})=\bigwedge
  ^pM$ for $c=2$ and {\it any} $p$.

  (2) For $c\ge 3$ the proof of Theorem \ref{BE's_guess} does not generalize
  since it is no longer true that
  $K_{{R/I_{t+1-p}(\varphi)}}\cong \bigwedge ^pM (\delta)$ for a suitable
  $\delta $ as the following example computed using Macaulay2 (\cite{Mac2})
  shows.
\end{remark}

\begin{example} \rm Set $R={\mathbb Q}[x_0,x_1,\cdots ,x_{n}]$ with $n$ either
  $13$ or $14$. Let $\cA$ be a $3\times 5$ matrix with general enough entries
  of linear forms. Take $M=\coker (\varphi ^* )$ where
  $\varphi ^*: R(-1)^5 \longrightarrow R^3$ is the graded morphism associated
  to $\cA$ and consider $I:=I_2(\cA)$ the ideal defined by the $2\times 2$
  minors of $\cA$. Using Macaulay2 we get that the minimal graded free
  $R$-resolutions of $R/I$, $K_{R/I}$ and $\bigwedge ^2 M$ have the following
  shape:
$$0 \rightarrow R(-10)^6\rightarrow R(-9)^{40}\rightarrow R(-8)^{105}\rightarrow \begin{array}{ccccccc} R(-6)^{50} \\ \oplus \\ R(-7)^{120} \end{array} \rightarrow \begin{array}{ccccccc} R(-5)^{168} \\ \oplus \\ R(-6)^{50} \end{array} \rightarrow R(-4)^{210} $$ $$ \rightarrow R(-3)^{120}\rightarrow R(-2)^{30}\rightarrow R\rightarrow R/I\rightarrow 0,$$

$$0 \rightarrow R(-11)\rightarrow R(-9)^{30}\rightarrow  R(-8)^{120} \rightarrow R(-7)^{210}
 \rightarrow \begin{array}{ccccccc} R(-5)^{50} \\ \oplus \\  R(-6)^{168}\end{array}
 \rightarrow \begin{array}{ccccccc}  R(-4)^{120} \\ \oplus \\ R(-5)^{50} \end{array} \rightarrow R(-3)^{105} $$ $$\rightarrow R(-2)^{40}\rightarrow R(-1)^6\rightarrow K_{R/I}\rightarrow 0,$$
and
$$0 \rightarrow R(-10)^3\rightarrow R(-9)^{15}\rightarrow  \begin{array}{ccccccc} R(-7)^{60} \\ \oplus \\ R(-8)^{15} \end{array} \rightarrow R(-6)^{210}
\rightarrow  R(-5)^{294} \rightarrow
 R(-4)^{210} \rightarrow  \begin{array}{ccccccc} R(-2)^{15} \\ \oplus \\ R(-3)^{60} \end{array} $$ $$\rightarrow R(-1)^{15}\rightarrow R^3\rightarrow \bigwedge^2M\rightarrow 0.$$

\noindent Therefore, $K_{R/I}$ is not isomorphic to a twist of $\bigwedge
^2M$. Nevertheless, using again Macaulay2 we compute a graded minimal $R$-resolution of $S_2(\bigwedge ^2M)$ and we obtain:
$$0 \rightarrow R(-10)\rightarrow R(-8)^{30}\rightarrow  R(-7)^{120} \rightarrow R(-6)^{210}
 \rightarrow \begin{array}{ccccccc} R(-4)^{50} \\ \oplus \\ R(-5)^{168} \end{array}
 \rightarrow
\begin{array}{ccccccc}
  R(-3)^{120}\\ \oplus \\ R(-4)^{50}
\end{array}
\rightarrow  R(-2)^{105} $$ $$\rightarrow R(-1)^{40}\rightarrow R^6\rightarrow S_2(\bigwedge^2M)\rightarrow 0.$$
We also have $K_{R/I}\cong S_2(\bigwedge^2M)(-1)$ (see Remark 3.11(3)).
\end{example}

The above example shows that, for $c\ge 3$,  it is no longer true that $K_{{R/I_{t+1-p}(\varphi)}}\cong \bigwedge ^pM (\delta)$  for a suitable $\delta $. Nevertheless, given $\varphi ^*:G^*\longrightarrow F^*$  a graded morphism between free $R$-modules of rank $t+c-1$ and $t$, respectively, there is a close relationship between the canonical module of $I_{t+1-p}(\varphi)$ and the Schur power $\Sigma ^IM$ of $M$ for a suitable partition $I$. Indeed, we have

\begin{theorem} \label{main2} Set $R=k[x_0,\cdots ,x_n]$. Let
  $\varphi ^*:G^*\longrightarrow F^*$ be a graded minimal morphism between free
  $R$-modules of rank $t+c-1$ and $t$, respectively. Assume that
  $\depth_{I_{t+1-p}(\varphi)} R = p(p+c-1)$. Set
  $\ell :=\sum _{j=1}^{t+c-1}a_j-\sum _{i=1}^tb_i$, $M=\coker (\varphi ^*)$
  and $J=(\overbrace{c-1,\cdots , c-1}^p)$. Then, we have:
$$K_{B_{t+1-p}}(n+1-p\ell) \cong \Sigma ^J M$$
\end{theorem}

\begin{proof} By hypothesis $M=\coker(\varphi ^*)$ and it has a presentation
$$G^*  \stackrel{\varphi ^*}{\longrightarrow }  F^*\longrightarrow M\longrightarrow 0.$$
Therefore, by \cite[Proposition V.2.2]{ABW}, $\Sigma ^JM$ has a presentation of the following type:
$$\Sigma ^{J_1}F^*\otimes G^*  \stackrel{\Sigma ^J \varphi ^*}{\longrightarrow } \Sigma ^J F^*\longrightarrow \Sigma ^J M\longrightarrow 0$$
where $J_1=(c-2,\overbrace{c-1,\cdots,c-1}^{p-1})$.

On the other hand, we know that the $\LL^{\bullet}(t+1-p,\varphi)$ gives us a minimal free $R$-resolution of $B_{t+1-p}$ of length $p(c+p-1)$ and its $R$-dual $\LL^{\bullet}(t+1-p,\varphi)^* $ provides a minimal free $R$-resolution of $K_{B_{t+1-p}}(n+1)$. Therefore, it is enough to describe the end of
$\LL^{\bullet}(t+1-p,\varphi)$.
 $\LL^{-p(c+p-1)}(t+1-p,\varphi)$ corresponds to the unique partition $I=(\overbrace{t+c-1, \cdots , t+c-1}^p)$ for which the Durfee square $p_I=p$, $n(I)=p(t-p)$ and $I'=(\overbrace{p,\cdots ,p}^{t-p},\overbrace{p+c-1,\cdots ,p+c-1}^p)$. Hence, we have $$\begin{array}{rcl} \LL^{-p(p+c-1)}(t+1-p,\varphi) & := & \Sigma ^{(I)^*}G^*\otimes \Sigma ^{I'}F \\
 & \cong & (\bigwedge ^{t+c-1} G^*)^{\otimes p}\otimes (\bigwedge ^tF)^{\otimes p}\otimes \Sigma ^ J F \\
 & \cong & \Sigma ^J F\otimes R(-p\ell ).\end{array}$$

\noindent Analogously, associated to
$\LL^{-p(p+c-1)+1}(t+1-p,\varphi)$ there is a unique partition, namely,  $I=(t+c-2,\overbrace{t+c-1, \cdots , t+c-1}^{p-1})$ with $p_I=p$, $n(I)=p(t-p)$ and $I'=(\overbrace{p,\cdots ,p}^{t-p},p+c-2,\overbrace{p+c-1,\cdots ,p+c-1}^{p-1})$. So,, we have $$\begin{array}{rcl} \LL^{-p(p+c-1)+1}(t+1-p,\varphi)& := & \Sigma ^{(I)^*}G^*\otimes \Sigma ^{I'}F \\ & \cong & (\bigwedge ^{t+1} G^*)^{\otimes (p-1)}\otimes \bigwedge ^t(G^*)\otimes (\bigwedge ^tF)^{\otimes p}\otimes \Sigma ^{J_1}  F \\
& \cong & G\otimes \Sigma ^{J_1} F\otimes R(-p\ell ) \end{array}$$
and we easily conclude arguing as in the last part of the proof of Proposition \ref{canonical module}.
\end{proof}

\begin{remark} \rm
(1) For $p=1$, we recover  Proposition \ref{resol} (iii), i.e. $K_{B_t}(n+1-\ell)\cong S_{c-1}M$.

(2) For $c=2$, we recover Proposition \ref{canonical module}, i.e. $K_{B_{t+1-p}}(n+1-p\ell)\cong \bigwedge ^p M$.

(3) For $t=3$ and $c=3$ we have $S_{c-1}(\bigwedge ^2M) \cong K_{B_{2}}(n+1-2\ell)$ where the last isomorphism follows from Theorem \ref{main2} and the Plethysm formula $$S_m(\bigwedge ^2F)=\bigoplus _{|I|=2m, i_j \text{ even for all j}} \Sigma ^I F.$$
\end{remark}

\section{The resolution of the exterior powers}

In Theorem \ref{BE's_guess}, we have found a minimal graded $R$-resolution of
$\bigwedge ^pM $ where $M=\coker \varphi^*$ and
$\varphi^*:G^*\longrightarrow F^*$ a graded minimal morphism between free
$R$-modules of rank $t-c+1$ and $t$, respectively, provided $c=2$. We already
pointed out that our approach does not generalize to $c\ge 3$. In this
section, supposing $\depth _J(B_t)\ge 2$, $J:=I_{t-1}(\varphi)$ (considered as
an ideal of $B_t$, i.e. we allow writing $J$ in place of $JB_t$), we will
outline an alternative proof that will provide a (possibly non-minimal) free
$R$-resolution of $\wedge ^2 M$ for $c=2, 3$ and an affirmative answer for
$c=p=2$, and for $p=2$ and $c=3$ up to a concrete summand which we guess
splits off (In fact, many cases have been checked with Macaulay2) of
Conjecture \ref{conjBE}. This approach applies several times the mapping cone
construction to the exact sequence
$$0\longrightarrow S_2M\otimes I_t\longrightarrow G^*\otimes M\stackrel{\varphi ^* \otimes id_M
}{\longrightarrow} F^*\otimes M \longrightarrow M\otimes M\longrightarrow 0
,$$
where $I_t:=I_{t}(\varphi)$ and then deduces a free resolution of
$\wedge ^2 M$ using the split-exact sequence
$$0\longrightarrow S_2M\longrightarrow M\otimes M\longrightarrow \bigwedge ^2
M\longrightarrow 0$$
and the free resolution of $S_2M$ given in Proposition \ref{resol}(ii). The
problem we encounter is to prove that
$\ker (\varphi ^* \otimes id_M) \cong S_2M\otimes I_t$ where we succeed for
$c=2,3$ assuming $\depth _J(B_t)\ge c+2$, but the method has the potential to
prove Conjecture \ref{conjBE} for $p=2$ and every $c \ge 2$ if we can show
$\depth _J(S_2M\otimes I_t)\ge 2$. Indeed if $c \ge 2$  and $\depth _J(B_t)\ge 2$ it is not so difficult is to see that
\begin{equation} \label{new1}
\ker (\varphi ^* \otimes id_M) \cong  \Hl _*^0(U,\widetilde{S_2M\otimes
  I_t}) \cong \Tor _1^R(M,M),  
\ \ \ \ {\rm where} \ \  U=\Proj(B_t)\setminus V(J).
\end{equation}
Since Macaulay2 computations indicate that
$\depth _J(S_2M\otimes I_t) \ge 2$ if $c \ge 2$ and $\depth _JB_t\ge c+2$, whence
$ S_2M\otimes I_t\ \cong \Hl _*^0(U,\widetilde{S_2M\otimes I_t})$, the method
may apply to any $c \ge 2$ with $B_t/J$ determinantal. To have sufficient
depth of  $S_2M\otimes I_t$ is related
to problems investigated in \cite{KM2015}, and we finish the paper by setting
forth another Conjecture \ref{newconj3}, seemingly unrelated to Conjecture
\ref{conjBE}, that will imply the above isomorphism.

Let us prove \eqref{new1}. Recall
$\ell_c:=\sum _{j=1}^{t+c-1}a_j-\sum _{i=1}^tb_i$ and that we  have the exact
sequence
\begin{equation}\label{exactseq2} G^* \stackrel{\varphi ^*}{\longrightarrow} F^* \longrightarrow M\longrightarrow 0. \end{equation}

\begin{remark}\label{continuation} \rm We can continue 
  (\ref{exactseq2}) by $P \stackrel{\epsilon}{\longrightarrow}
  G^*\longrightarrow \cdots $ to the left where $P:=\bigwedge ^{t+1}G^*\otimes
  \bigwedge
  ^tF$ and $\epsilon $ is the "splice map" (see Proposition \ref{resol}(ii)).
  Let $-1 \le i \le c$. Since $\epsilon \otimes _R id_{S_iM}=0$ and  $ \Tor
  _1^R(M,S_iM) := \ker (\varphi ^* \otimes id_{S_iM}) / \im(\epsilon\otimes
  _Rid_{S_iM})$, we get that the sequence
$$ 0\longrightarrow \Tor _1^R(M,S_iM) \longrightarrow G^* \otimes S_iM
\stackrel{\varphi ^* \otimes id_{S_iM}}{\longrightarrow} F^*\otimes S_iM
\longrightarrow M\otimes S_iM\longrightarrow 0 $$
is exact. Thus $i=1$ takes care of one isomorphism in \eqref{new1}. Moreover
note that if $\depth _J(B_t)\ge 2$, the above exact sequence shows
$\depth _J \Tor^R_{1}(M,S_iM) \ge 2$ as  $S_iM$ is maximally Cohen-Macaulay.
\end{remark}

\begin{lemma} \label{S2M} Let $U=X\setminus V(J)$ where $J=I_{t-1}(\varphi)$, $X=\Proj(B_t)$ and suppose $B_t$ is standard determinantal and that $\depth _J(B_t)\ge 2$. It holds
$$\widetilde{\Tor _1^R(M,M)}\cong \widetilde{ S_2M\otimes _RI_t} \text{ as sheaves over } U.$$

\noindent In particular,
$\Tor _1^R(M,M)\cong \Hl _*^0(U,\widetilde{S_2M\otimes I_t})$.
\end{lemma}

\begin{proof} First note that tensoring $0 \to I_t \to R \to B_t \to 0$ with
  the $B_t$-module $M$ and using the induced long exact sequence of $\Tor_j$
  for $j \le 1$ we easily get $\Tor^R_{1}(M,B_t) \cong M\otimes _RI_t$. Then
  consider the sheafification over $U$ of the exact sequence in Remark
  \ref{continuation} for $i=0$. Since $\widetilde{M_{|U}}$ is invertible,
  hence flat over $U$, we can tensor this sequence with $\widetilde{M_{|U}}$,
  and we get the sheafification of the exact sequence of Remark
  \ref{continuation} for $i=1$. This shows the first isomorphism in
$$\widetilde{\Tor_1^R(M,M)}\cong \widetilde{M} \otimes _{\tilde B_t}
\widetilde{\Tor_1^R(M,B_t)} \cong \widetilde{M} \otimes _{\tilde B_t}
(\widetilde{M\otimes _R I_t)}\cong \widetilde{S_2M\otimes _R I_t}.$$
Since $\depth _J\Tor_1^R(M,M)\ge 2$ by Remark \ref{continuation}, we also
get the final conclusion.
\end{proof}

Combining Remark \ref{continuation} and Lemma \ref{S2M}, we get that
\eqref{new1} holds.

\begin{lemma}\label{lem3} Suppose $B_t$ is standard determinantal.

(i)  If $\depth _J B_t\ge 4$ and $c=\codim _R B_t=2$, then
$$\Tor_1^R(M,M)\cong  S_2M\otimes _RI_t.$$

 (ii) If $\depth _J B_t\ge 5$ and $c=\codim _R B_t=3$, then
$$\Tor_1^R(M,M)\cong  S_2M\otimes _RI_t \cong \Ext^2_R(I_t,I_t)(-\ell_3) \quad \text{ where } \quad \ell_3=\sum
_{j=1}^{t+2}a_j-\sum _{i=1}^tb_i .$$

\end{lemma}

\begin{proof}
  (i) As $c = 2$ and $\depth_J B_t \ge 4$, \cite[Corollary 3.7]{KM2015}
  applies and we get $$\Ext_R^1(B_t,M) \cong \Ext_R^1(M,S_2M)$$ by choosing
  $r = 1$, $s =2$ and $a=1$ in Corollary 3.7. Using \eqref{torext} below with
  $j = 1$ and $(r,s)=(0,2)$, and also with $(r,s)=(1,1)$, we get
  $\Tor_1^R(M,M) \cong \Tor_1^R(B_t,S_2M)$. But if we tensor
  $0 \to I_t \to R \to B_t \to 0$ with the $B_t$-module $S_2M$ we deduce
  $\Tor^R_{1}(B_t,S_2M) \cong S_2M\otimes _RI_t$. And combining the two last
  isomorphisms we conclude as claimed.

(ii) In this case the first isomorphism in (ii) follows from \cite[Corollary
3.7]{KM2015}. Indeed choosing $r = s = 1$ and $a=2$ in the
``resp.-assumptions'' of Corollary 3.7, we get
$$\Ext_R^2(M,M) \cong \Ext_R^2(B_t,B_t).$$ But \eqref{torext} below with
$j = 1$ and $r = s = 1$ yields $\Tor_1(M,M)(l) \cong \Ext^2(M,M)$ where $\ell=\ell_3$, and
\eqref{torext} with $j=1$, $r = 2$ and $s = 0$ yields
$\Ext^2(B_t,B_t) \cong \Tor_1(S_2M,B_t)(l)$. Thus we get
$$\Tor_1^R(M,M) \cong \Tor_1^R(S_2M,B_t)$$ which proves what we want since
$ \Tor_1(S_2M,B_t) \cong S_2M \otimes I_t$ can be shown as in (i).

As $c=3$, we have $\pd I_t=2$, whence $\Ext^2_R(I_t,-)$ right-exact and since
$S_2M$ is by Proposition \ref{resol}(iii) the canonical module
$K_{B_t} = \Ext^2_R(I_t,R)(-n-1)$ of $B_t$ (up to twist), we get {\small
$$\begin{array}{rcc}
    0=\Ext^1_R(I_t,R)\longrightarrow\Ext^1_R(I_t,B_t)\longrightarrow\Ext_R^2(I_t,I_t)\longrightarrow\Ext^2_R(I_t,R) & \longrightarrow & \Ext^2_R(I_t,B_t) \\
    \downarrow \cong  & & \downarrow \cong \\
    S_2M(\ell) & \cong &  S_2M(\ell)\otimes B_t .  \end{array}$$
}Indeed to see the twist, if we apply $\Hom_R(-,R)$ to the $R$-free resolution
$(P_{\bullet })$ of $I_t$ in
\eqref{dualS2} below, we find exactly the $R$-free resolution of $S_2M(\ell)$
given by \eqref{splice} and
Proposition \ref{resol}(ii). Moreover as the right-exactness of  $\Ext^2_R(I_t,-)$ also implies  $\Ext^2_R(I_t,I_t)
\cong  \Ext^2_R(I_t,R) \otimes I_t$, we get
\begin{equation} \label{S2MotimesI}
\Ext^1_R(I_t,B_t)\cong \Ext^2_R(I_t,I_t)\cong S_2M(\ell_3)\otimes I_t,
\end{equation}
as desired.
\end{proof}

\begin{remark} \label{IBt} \rm For generic determinantal rings (i.e. where the
  matrix $\cA$ consists of the indeterminites of the polynomial ring R) we
  have used Macaulay2 over the field of rational numbers to show
  $\pd (S_2M\otimes _RI_t) = 2c$ for $2 \le c \le 5$ and $t = 3$ as well as
  for $2 \le c \le 9$ and $t = 2$. Using notably Auslander-Buchsbaum formula
  connecting depth and projective dimension, we get
  $\depth_J (S_2M\otimes _RI_t) \ge 2$ as $\dim B_t = c+2$, hence
  $\dim R = 2c+2$. This indicates that Lemma \ref{lem3} likely extends to any
  $c \ge 2$ as Conjecture~\ref{newconj3} asserts. So these Macaulay2
  computations give some evidence to an important case of that conjecture,
  namely that $\Tor^R_{1}(S_2M,B_t) \cong \Tor^R_{1}(M,M)$ holds. As
  $\Tor^R_{1}(S_2M,B_t) \cong S_2M\otimes _RI_t$, Conjecture~\ref{newconj3}
  asserts
\begin{equation} \label{Tor1S2}
\Tor^R_{1}(M,M) \cong S_2M\otimes _RI_t
\end{equation}
provided   $\depth _JB_t \ge c+2$ which by Remark \ref{continuation} implies
the exactness of the desired \eqref{new1}. Note that only assuming $\depth
_JB_t \ge 2$ and $B_t$ determinantal, \eqref{Tor1S2} is
equivalent to $\depth _J(S_2M\otimes _RI_t) \ge 2$ by Lemma \ref{S2M} as
  $\depth _J \Tor^R_{1}(M,M) \ge 2$ by Remark \ref{continuation}. Thus
  $\pd(S_2M\otimes _RI_t) = 2c$ implies  \eqref{Tor1S2}.
\end{remark}

As a consequence of the above results we will be able to prove Conjecture
\ref{conjBE} for $c=p=2$, and for $p=2$ and $c=3$ up to a concrete direct summand.
 As a main tool we will use
the following generalization of the mapping cone process:

\begin{lemma}\label{technicallemma}
  Let
  $Q_{\bullet} \stackrel{\sigma _{\bullet}}{\longrightarrow} P _{\bullet}
  \stackrel{\tau _{\bullet}}{ \longrightarrow} F _{\bullet} $
  be morphisms of complexes satisfying $Q_j=P_j=F_j=0$ for $j < 0$. Denote by
  $d_Q^{i}:Q_i \to Q_{i-1}$ (respectively $d_P^{i}$, $d_F^{i}$) the
  differentials of $Q_{\bullet}$ (respectively $P_{\bullet}$,
  $F_{\bullet}$). Assume that all three complexes are acyclic (exact for
  $j \ne 0$) and that the sequence
\begin{equation*}
  0 \longrightarrow \coker d_Q^{1}  \longrightarrow \coker d_P^{1} \stackrel {\alpha}
  { \longrightarrow}  \coker d_F^{1}
\end{equation*}
is exact. Moreover assume that there exists a morphism $\ell _{\bullet}: Q_{\bullet}\longrightarrow  F_{\bullet}[1]$ such that
for any integer $i$:
\begin{equation} \label{ellcomp} d_F^{i+1}\ell _i+\ell _{i-1}d_Q^{i}=\tau
  _i\sigma_i.
\end{equation}
Then, the complex  $Q_{\bullet}\oplus P_{\bullet}[1]\oplus F_{\bullet}[2]$ given by
$$Q_{i}\oplus P_{i+1}\oplus F_{i+2} \stackrel {d^{i}_{Q,P,F}}
{ \longrightarrow} Q_{i-1}\oplus P_{i}\oplus F_{i+1}$$ where
 $$d^{i}_{Q,P,F}:=\begin{pmatrix}d^{i}_Q & 0 & 0 \\ \sigma _{i} & -d_P^{i+1} & 0 \\
 \ell _{i} & -\tau_{i+1} & d^{i+2}_F \end{pmatrix} $$ is acyclic (exact for $i \ne -2$) and $ \coker  {d^{-1}_{Q,P,F}} = \coker \alpha$.
 \end{lemma}

\begin{proof} See \cite[ Lemma 3.1]{KM2015}.
\end{proof}

We are now ready to prove the main result of this section.

\begin{theorem} \label{resolwedge2} Let $J=I_{t-1}(\varphi)$ and  $B_t$ be standard determinantal, and assume either
  $\depth _JB_t\ge 2$ and $\depth _J(S_2M\otimes _RI_t) \ge 2$, or
   $\depth _JB_t\ge c+2$ and $2\le c\le 3$.
  Then Conjecture
  \ref{conjBE} (with $\varphi $ minimal) holds for $c=p=2$, and for $p=2$ and $c\ge 3$ up to a direct summand.
\end{theorem}

\begin{proof} As we have $\Tor_1^R(M,M)\cong S_2M\otimes _RI_t$ by Lemma
  \ref{S2M}, Lemma \ref{lem3} and assumptions, the sequence
\begin{equation} \label{mains2m} 0\longrightarrow S_2M\otimes I_t\longrightarrow
  G^*\otimes M\longrightarrow F^*\otimes M \longrightarrow M\otimes
  M\longrightarrow 0
\end{equation}
is exact.
Using the free resolution of $S_2M$ given in Proposition \ref{resol}(ii)
\begin{equation} \label{s2m} \cdots \longrightarrow \bigwedge ^{t+2} G^*\otimes \bigwedge ^t
F\longrightarrow \bigwedge ^2G^* \longrightarrow G^*\otimes F^*\longrightarrow
S_2F^*\longrightarrow S_2M\longrightarrow 0
\end{equation}
and the Eagon-Northcott resolution of $I_t$ appearing in the horizontal row
below
$$
 \begin{array}{cccccccccc}
    &  & \bigwedge ^tG^*\otimes \bigwedge ^tF\otimes F^*\otimes  G^* &&& \\
   \downarrow && \downarrow &&& \\
  \bigwedge ^{t+1}G^*\otimes F \otimes \bigwedge ^tF\otimes S_2F^* & \longrightarrow & \bigwedge ^tG^* \otimes \bigwedge ^tF\otimes S_2F^* &&& \\
   \downarrow && \downarrow &&& \\
 (\bigwedge ^{t+1}G^*\otimes F \otimes \bigwedge ^t F)\otimes S_2M & \longrightarrow & (\bigwedge ^tG^* \otimes \bigwedge ^tF) \otimes S_2M & \longrightarrow & S_2M\otimes I_t & \longrightarrow 0 \\
 \end{array}
$$
it is easy to see that a (possibly non-minimal) presentation of
$S_2M\otimes _RI_t$  starts out as:
$$
(G^*\otimes F^*\otimes \bigwedge ^tG^*\otimes \bigwedge^tF)\oplus (S_2F^*\otimes \bigwedge ^{t+1}G^* \otimes F\otimes \bigwedge ^t F)  \longrightarrow  S_2F^*\otimes \bigwedge ^tG^*\otimes \bigwedge ^tF  \longrightarrow S_2M\otimes I_t \longrightarrow 0.
$$
Free  $R$-resolutions of $ G^*\otimes M$ and of $F^*\otimes M$
are given by Proposition \ref{resol}(ii). Therefore, applying the generalization of the mapping cone process (Lemma \ref{technicallemma}) we obtain that a free $R$-resolution of $M\otimes M$ begins as follows:
{\small $$\begin{array}{ccccccc} F^*\otimes \bigwedge ^tF\otimes F\otimes \bigwedge ^{t+2}G^* \\ \oplus \\ G^*\otimes \bigwedge ^tF\otimes \bigwedge ^{t+1} G^* \\ \oplus \\
G^*\otimes F^*\otimes \bigwedge ^tG^*\otimes \bigwedge^tF\\ \oplus  \\S_2F^*\otimes \bigwedge ^{t+1}G^* \otimes F\otimes \bigwedge ^t F
\end{array}\longrightarrow
\begin{array}{ccccccc} G^*\otimes G^* \\ \oplus \\
F^*\otimes \bigwedge ^tF\otimes \bigwedge ^{t+1} G^* \\ \oplus \\ S_2F^*
  \otimes \bigwedge ^tG^* \otimes  \bigwedge ^tF \end{array}\longrightarrow
(F^*\otimes G^*)^{\oplus 2}\longrightarrow F^*\otimes F^* \longrightarrow M\otimes M\longrightarrow 0.
$$}

Using Pieri's formula  \ref{Pieri2} and letting $I=(0,\overbrace{1,\cdots ,1}^{t-2},3)$ we have
{\small $$\begin{array}{rcl} S_2F^*\otimes \bigwedge ^{t+1}G^* \otimes F\otimes \bigwedge ^t F & \cong & (S_2F^*\otimes \bigwedge ^{t-1}F^*)\otimes \wedge ^tF\otimes \bigwedge ^{t+1}G^* \otimes \bigwedge ^t F \\ & \cong &
(F^*\otimes \bigwedge ^tF^*\oplus \Sigma ^{I}F^*)\otimes \bigwedge ^tF\otimes \bigwedge ^{t+1}G^* \otimes \bigwedge ^t F \\
 & \cong & F^*\otimes \bigwedge ^tF\otimes \bigwedge ^{t+1} G^*  \oplus \Sigma ^{I}F^*\otimes \bigwedge ^tF\otimes \bigwedge ^{t+1}G^* \otimes \bigwedge ^t F.
\end{array}
$$}
It is important to point out that the  summand
$H:=F^*\otimes \bigwedge ^tF\otimes \bigwedge ^{t+1} G^*$ will split off if the corresponding matrix representing the last morphism of the above free $R$-resolution of $M\otimes M$ is invertible. Macaulay2 computations and explicit description of the involved morphism for small 
$t$ give evidences that it should be true.

Using \eqref{s2m},  the short exact sequence
$$0\longrightarrow S_2M\longrightarrow M\otimes M\longrightarrow \bigwedge ^2 M\longrightarrow 0,$$
and applying once more the mapping cone process we get that a free $R$-resolution of $\bigwedge ^2 M$ looks like

{\small $$\begin{array}{ccccccc} F^*\otimes \bigwedge ^tF\otimes F\otimes \bigwedge ^{t+2}G^* \\ \oplus \\ G^*\otimes \bigwedge ^tF\otimes \bigwedge ^{t+1} G^* \\ \oplus \\
G^*\otimes F^*\otimes \bigwedge ^tG^*\otimes \bigwedge^tF\\ \oplus  \\ \Sigma ^{I}F^*\otimes \bigwedge ^tF\otimes \bigwedge ^{t+1}G^* \otimes \bigwedge ^t F \\  \oplus \\ H
\end{array}\longrightarrow
\begin{array}{ccccccc} S_2  G^*  \\ \oplus \\ S_2F^* \otimes \bigwedge ^tG^* \otimes  \bigwedge ^tF \\  \oplus \\ H \end{array}\longrightarrow F^*\otimes G^*\longrightarrow  \bigwedge ^2 F^* \longrightarrow \bigwedge ^2 M\longrightarrow 0
$$}
and hence, except for $H$, it starts as predicted in Conjecture \ref{conjBE}; cf. Theorem 3.7 for the case $c=p=2$.
\end{proof}

\vskip 2mm The last part of this section has as a main goal to explicitly
construct a full free $R$-resolution (possibly non-minimal) of $\bigwedge ^2M$
for $2\le c\le 3$ using mapping cone procedures. As an intermediate result we
obtain a free graded resolution of the normal
module $N_{B_t}$ for $c=3$ (the case $c=2$ is well known). A problem with the
above procedure is that the free resolution often contains redundant summands.
In some cases one may more easily get rid of all redundant summands using a generalization
of the mapping cone process. We illustrate this in Theorem 4.9 by
finding a minimal free graded resolution of the normal module $N_{B_t}$ for
$c=3$. We analyze separately 2 cases:

\vskip 2mm
\noindent \underline{Case 1:} Suppose $c=2$ and $\depth _J B_t\ge 4$. Applying
$\Hom_R(-,S_2M)$ to the exact sequence (\ref{exactseq2}) continued by
$P \stackrel{\epsilon}{\longrightarrow} G^*\longrightarrow \cdots $  to the
left as in Remark \ref{continuation}, $\epsilon $ the "splice
map", and using $\Ext^1_R(M,S_2M)\cong S_2M\otimes I_t(\ell_2)$ (which follows from
 \eqref{torext} below) we get the following exact sequence
\begin{equation} \label{crellec2}
0\longrightarrow  M  \longrightarrow  F\otimes S_2M  \longrightarrow  G\otimes S_2M  \longrightarrow   S_2M\otimes
I_t(\ell_2)  \longrightarrow 0
\end{equation}
as $\Hom_R(\epsilon,S_2M ) = 0$, (e.g. see \cite[Equation (3.6)]{KM2015}  for
details). Since $M$ and $S_2M$ are maximally Cohen-Macaulay, the mapping cone
construction applied to (\ref{crellec2}) yields a free $R$-resolution of
$S_2M\otimes I_t(\ell_2)$. Using the exact sequence (\ref{mains2m}), this leads via mapping cone to a free $R$-resolution of
$M\otimes M$ and using the short exact sequence
$0\longrightarrow S_2M\longrightarrow M\otimes M\longrightarrow \bigwedge
^2M\longrightarrow 0$
we finally get a free $R$-resolution of $\bigwedge ^2M$ of length $6$; 
cf. Theorem \ref{BE's_guess} which in fact proves
more.

\vskip 2mm
\noindent \underline{Case 2:} Suppose $c=3$ and $\depth _J B_t\ge 5$. Recall
that the Eagon-Northcott resolution of $I_t$ yields
\begin{equation} \label{dualS2}
 (P_{\bullet }): \quad 0\longrightarrow S_2F(-\ell)\longrightarrow G\otimes
 F(-\ell)\longrightarrow \bigwedge ^2G(-\ell)\longrightarrow
 I_t\longrightarrow 0
\end{equation}
 with
$\ell:=\ell_3 =\sum _{j=1}^{t+2}a_j-\sum _{i=1}^tb_i$, as $ \bigwedge
^{t+i}G \otimes \bigwedge ^tF \cong  \bigwedge ^{2-i}G(-\ell)\ $. 
Firstly, we describe a
minimal free $R$-resolution of the normal module
$N_{B_t}:=\Hom_R(I_t,B_t)\cong \Ext^1_R(I_t,I_t)$. As a main tool we will use again the generalization of the mapping cone process. We have:

\begin{theorem}\label{resolminimlanormal}
Let $\varphi: F:=\oplus _{i=1}^tR(b_i)\longrightarrow G:=\oplus _{j=1}^{t+2} R(a_j)$ be a graded minimal $R$-morphism between two free $R$-modules of rank $t$ and $t+2$, respectively. Set  $\ell:=\ell_3=\sum_{j=1}^{t+2}a_j-\sum_{i=1}^tb_{i}$ and assume $\depth _JB_t\ge 4$ where $J=I_{t-1}(\varphi )$. Then the normal module
 $N_{B_t}:=\Hom_R(I_t,B_t)\cong \Ext^1_R(I_t,I_t)$ has a minimal free $R$-resolution of the following type:
 $$ (N_{\bullet}): \quad
0 \longrightarrow \bigwedge ^2F(-\ell) \longrightarrow  (G\otimes F(-\ell))\longrightarrow (F\otimes G^*)\oplus S_2G(-\ell)$$
$$\longrightarrow((F\otimes F^*)\oplus (G\otimes G^*))/R\longrightarrow G\otimes F^*\longrightarrow N_{B_t}\longrightarrow 0.
$$
\end{theorem}
\begin{proof}
 Since $\depth _JB_t\ge 4$, using \cite[Corollary 4.9 and Theorem 5.2]{K2014}, we get $N_{B_t}\cong \Ext^1_R(M,M)$. Then, \cite[Theorem 3.1 ]{K2014} implies $\depth N_{B_t}=\depth \Ext^1_R(M,M)\ge \dim B_t-1$ and its proof implies that the sequence
\begin{equation}\label{normalmodule}0\longrightarrow B_t\longrightarrow F\otimes M \longrightarrow G\otimes M\longrightarrow N_{B_t}\longrightarrow 0\end{equation}
is exact. This last exact sequence together with Lemma \ref{technicallemma} will allow us to construct a minimal free $R$-resolution of the normal module $N_{B_t}$. Indeed,
the Buchsbaum-Rim resolution (see Proposition \ref{resol}):
$$ 0\longrightarrow F(-\ell)\cong \bigwedge ^{t+2}G^*\otimes F\otimes \bigwedge ^tF\longrightarrow \bigwedge ^{t+1}G^*\otimes  \bigwedge ^tF\cong G(-\ell)\longrightarrow G^*\longrightarrow F^* \longrightarrow M \longrightarrow 0$$
together with the resolution $ (P_{\bullet })$ of $B_t$ yields the diagram

\[
\xymatrix{ &0 \ar[d]  &0 \ar[d] &0 \ar[d]  &\\
&  S_2F(-\ell)\ar[d]  \ar[r]_{\tau _2} & F\otimes F(-\ell)\ar[d] \ar[r] & G\otimes F(-\ell)\ar[d] & \\
& G\otimes F(-\ell)\ar[d] \ar[r]_{ \tau _1} \ar[rru]^{\quad \quad \ell _1} & F\otimes G(-\ell)\ar[d]  \ar[r] &  G\otimes G(-\ell) \ar[d] & \\
  & \bigwedge ^2G(-\ell)\ar[d] \ar[r]_{\tau _0} \ar[rru]^{\quad \quad \ell _0} & F\otimes G^* \ar[d]  \ar[r] & G\otimes G^* \ar[d] & \\
 & R  \ar[d] \ar[r]_{\tau _{-1}} \ar[rru]^{\quad \quad \ell _{-1}}& F\otimes F^* \ar[d]  \ar[r] & G\otimes F^* \ar[d] & \\
 0\ar[r] &  B_t \ar[d]  \ar[r] & F\otimes M \ar[d]  \ar[r] & G\otimes M \ar[d]   \ar[r] & N_{B_t}  \ar[r] &0 \\
&0 &0 &0 &}
\]
where $\ell _1$, $\ell _0$, $\ell _{-1}$, $\tau _2$, $\tau _1$ and $\tau _{-1}$ are the natural injections. Since $\coker (\tau _2)\cong \bigwedge ^2F(-\ell)$,  $\coker (\ell _0)\cong S_2G(-\ell)$  and $\tau _1=\id$, the mapping cone construction (see Lemma \ref{technicallemma})
gives us the minimal free $R$-resolution  of $N_{B_t}$ (it is minimal as all vertical and rightmost horizontal maps in the above diagram are minimal):
$$(N_{\bullet}): \quad
0 \longrightarrow \bigwedge ^2F(-\ell) \longrightarrow G\otimes F(-\ell)\longrightarrow (F\otimes G^*)\oplus S_2G(-\ell)$$
$$\longrightarrow((F\otimes F^*)\oplus (G\otimes G^*))/R\longrightarrow G\otimes F^*\longrightarrow N_{B_t}\longrightarrow 0
$$
which proves what we want.
\end{proof}

Let us now compute a free $R$-resolution of $S_2M\otimes I_t$. We split $(P_{\bullet })$ into two short exact sequences:
 \begin{equation}\label{P1}  0\longrightarrow P_{-3}:= S_2F(-\ell)\longrightarrow P_{-2}:= G\otimes F(-\ell)\longrightarrow E\longrightarrow 0\end{equation}
 \begin{equation}\label{P2} 0\longrightarrow E\longrightarrow P_{-1}:= \bigwedge ^2G(-\ell)\longrightarrow I_t\longrightarrow 0\end{equation}
 and we use Lemma \ref{lem3} which implies $S_2M(\ell )\otimes I_t\cong \Ext^2_R(I_t,I_t)$. By (\ref{P2}) we have $\Ext^2_R(I_t,I_t)\cong \Ext^1_R(E,I_t)$, and
 {\small
 \begin{equation}\label{diagram}
 \begin{array}{ccccccccccc}
   0\longrightarrow & \Hom(I_t,I_t) & \longrightarrow & \Hom(P_{-1},I_t) & \longrightarrow & \Hom(E,I_t) & \longrightarrow & \Ext^1_R(I_t,I_t) & \longrightarrow 0 \\
                    & \cong\uparrow && \cong \uparrow & & \parallel & & \cong\uparrow \\
                    & R & \longrightarrow & P_{-1}^*\otimes I_t & \longrightarrow & \Hom(E,I_t) & \longrightarrow & N_{B_t} & \longrightarrow 0 \\
                    & & & \uparrow & & & & \uparrow \\
                    & & &  P_{-1}^*\otimes P_{-1} & & & & G\otimes F^* \\
                    & & & \uparrow & & & & \uparrow \\
                    & & &  P_{-1}^*\otimes P_{-2} & & & & ((F\otimes F^*)\oplus (G\otimes G^*))/R \\
                    & & & \uparrow & & & & \uparrow \\
                    & & &  P_{-1}^*\otimes P_{-3} & & & & (F\otimes G^*)\oplus S_2G(-\ell ) \\
                    & & & \uparrow & & & & \uparrow \\
                    & & &  0 & & & & G\otimes F(-\ell) \\
                    & & & & & &  & \uparrow \\
                    & & &  & & &  & \bigwedge ^2F(-\ell) \\
                    & & & & & & & \uparrow \\
                    & & &  & & &  & 0
 \end{array}
 \end{equation}
} whence
 a free resolution of $\Hom(E,I_t)$ is:
$$ 0 \longrightarrow \bigwedge ^2F(-\ell)\longrightarrow G\otimes F(-\ell) \longrightarrow (F\otimes G^*)\oplus S_2 G(-\ell) \oplus (P_{-1}^*\otimes P_{-3}) \longrightarrow $$ $$ ((F\otimes F^*)\oplus (G\otimes G^*))/R \oplus (P_{-1}^*\otimes P_{-2}) \longrightarrow (G\otimes F^*)\oplus  (P_{-1}^*\otimes P_{-1}/R) \longrightarrow \Hom(E,I_t) \longrightarrow 0.$$
Now, we apply the functor $\Hom(-,I_t)$ to the exact sequence (\ref{P1}) and we get
{\small
 $$
 \begin{array}{ccccccccccc}
 0\longrightarrow  \Hom(E,I_t) & \longrightarrow & \Hom(P_{-2},I_t) & \longrightarrow & \Hom(P_{-3},I_t) & \longrightarrow & \Ext^1_R(E,I_t)  \longrightarrow 0. \\
  \parallel &&  \cong \uparrow & & \cong \uparrow & & \cong\uparrow \\
 0\longrightarrow \Hom(E,I_t) & \longrightarrow & P_{-2}^*\otimes I_t & \longrightarrow & P_{-3}^*\otimes I_t & \longrightarrow & \Ext^2_R(I_t,I_t) \\
    \uparrow & & \uparrow & &  \uparrow \\
    (G\otimes F^*)\oplus  (P_{-1}^*\otimes P_{-1}/R) &\longrightarrow &  P_{-2}^*\otimes P_{-1} & \longrightarrow & P_{-3}^*\otimes P_{-1} \\
     \uparrow & & \uparrow & & \uparrow \\
      ((F\otimes F^*)\oplus (G\otimes G^*))/R \oplus (P_{-1}^*\otimes P_{-2}) &\longrightarrow &  P_{-2}^*\otimes P_{-2} &\longrightarrow & P_{-3}^*\otimes P_{-2}  \\
    \uparrow & & \uparrow & &  \uparrow \\
     (F\otimes G^*)\oplus S_2G(-\ell) \oplus (P_{-1}^*\otimes P_{-3}) &\longrightarrow &  P_{-2}^*\otimes P_{-3} & \longrightarrow & P_{-3}^*\otimes P_{-3} \\
     \uparrow & & \uparrow & &  \uparrow & & &  \\
     G\otimes F(-\ell) & & 0  & & 0 & &  \\
     \uparrow & & & & &   \\
     \bigwedge ^2F(-\ell) & &  & & &  &  \\
     \uparrow & & & & & &  \\
     0 & &  & & &  &
 \end{array}
 $$
} Finally, applying the mapping cone process we obtain the following free $R$-resolution of $\Ext^2(I_t,I_t)(-\ell)\cong S_2M\otimes I_t$:
{\small
$$0\longrightarrow \bigwedge ^2F(-2\ell)\longrightarrow G\otimes F(-2\ell)\longrightarrow (F\otimes G^*)(-\ell)\oplus S_2 G(-2\ell) \oplus (P_{-1}^*\otimes P_{-3})(-\ell)\longrightarrow $$
$$((F\otimes F^*)\oplus (G\otimes G^*))/R \oplus (P_{-1}^*\otimes P_{-2})(-\ell)\oplus (P_{-2}^*\otimes P_{-3})(-\ell)\longrightarrow  $$
$$(G\otimes F^*)\oplus ((P_{-1}^*\otimes P_{-1})\oplus (P_{-2}^*\otimes P_{-2})\oplus (P_{-3}^*\otimes P_{-3}))/R(-\ell)\longrightarrow $$
$$\begin{array}{ccccc} (P_{-2}^*\otimes P_{-1})(-\ell)\oplus (P_{-3}^*\otimes P_{-2})(-\ell) & \longrightarrow & (P_{-3}^*\otimes P_{-1})(-\ell) & \longrightarrow S_2M\otimes I_t \longrightarrow 0.\\
    \parallel && \parallel \\
    (G^*\otimes F^*\otimes \bigwedge ^2 G)(-\ell) 
    \oplus (S_2F^*\otimes G \otimes F)(-\ell)
 &
   \longrightarrow & S_2F^*\otimes \bigwedge ^2 G(-\ell)
\end{array}$$
}Once we have a free $R$-resolution of $S_2M\otimes I_t$ arguing as in Case
1, see also the proof of  Theorem~\ref{resolwedge2},
we get a free $R$-resolution of $\bigwedge ^2M$ of length $8$.
\begin{remark}\rm Using Macaulay2 we have checked that the above free $R$-resolution of $S_2M\otimes I_t$ can contain redundant summands.
\end{remark}


\section{Final comments}
Based on the previous results and on a series of examples computed with Macaulay2, we state the following guesses:

\begin{conj} \label{newconj2} Set $R=k[x_0,\cdots ,x_n]$. Let
  $\varphi ^*:G^*\longrightarrow F^*$ be a graded morphism between free
  $R$-modules of rank $t+c-1$ and $t$, respectively. Assume that
  $\depth_{I_{t+1-p}} R=p(c+p-1)$ and $c \ge 2$. Then for
  any $1\le i \le c+p-1$, $S_i(\bigwedge ^p(M))$ is a maximal Cohen-Macaulay
  $R/I_{t+1-p}$-module.
\end{conj}

\begin{remark} \rm
(1)  By \cite[Proposition 4]{V}, we know that Conjecture \ref{newconj2} is true for $i=1$ and any $t$, $c$ and $p$.

(2) Using Macaulay2 we have checked Conjecture \ref{newconj2} for any $1\le i \le c+p-1$ and any $p$ provided $c\le 4$ and $t\le 5$ and found it to be true.
\end{remark}

\begin{conj} \label{newconj3} Set $R=k[x_0,\cdots ,x_n]$ and let
  $B_t:=R/I_{t}(\varphi)$ be standard determinantal with $\depth _JB_t\ge 2$,
  $J := I_{t-1}(\varphi)$ where $\varphi:F \to G$ is a graded morphism
  between free $R$-modules of ranks $t$ and $t+c-1$ respectively. Let
  $j, r, s$ and $q$ be integers and assume
  $j>\max(0,\frac{c+3-\depth_JB_t}{2})$ and $c \ge 2$. Then for
  $r, s, r-q,s+q \in [-1, c]$, it holds
  \begin{equation} \label{conj53eq}
  \Tor^R_{j}(S_{r}M,S_{s}M) \cong \Tor^R_{j}(S_{r-q}M,S_{s+q}M).
  \end{equation}
\end{conj}

\begin{remark} \label{tonewconj3} \rm Using Macaulay2 we have checked that the
  hypothesis $j>\max(0,\frac{c+3-\depth_JB_t}{2})$ in Conjecture
  \ref{newconj3} cannot be dropped, at least not when
  $\depth_JB_t = \dim B_t < c+2$. Indeed,

(1) We consider the polynomial ring $R=\QQ[a,b,c,d,e,f,g]$ and $\cA$ a
$3\times 5 $ matrix with entries general linear forms in $R$. Associated to
$\cA$ we have  a graded morphism  $\varphi $ between free $R$-modules of rank
$t=3$ and $t+c-1=5$. We have got  that $\Tor^R_{2}(M,M)$ and
$\Tor^R_{2}(S_{2}M,B_t)$ have the same graded Betti numbers and thus are likely
isomorphic. Moreover we got $\pd \Tor^R_{1}(S_{2}M,B_t) = 6$, thus $\depth
\Tor^R_{1}(S_2M,B_t)= 1$ by Auslander-Buchsbaum's formula. As  $\depth
\Tor^R_{1}(M,M)\ge 2$ by Remark \ref{continuation},  $\Tor^R_{1}(S_2M,B_t)$
cannot  be isomorphic to  $\Tor^R_{1}(M,M)$. Note  $\depth_JB_t = \dim B_t = 4$.

(2) Analogous result replacing $R$ by $\QQ[a,b,c,d,e,f]$, but now
$\pd \Tor^R_{1}(S_{2}M,B_t) =6$ yields $\depth \Tor^R_{1}(S_2M,B_t)= 0$ and
thus $\Tor^R_{1}(S_{2}M,B_t) \not\cong \Tor^R_{1}(M,M)$ while
$\depth_JB_t = \dim B_t = 3$.

(3) Taking $R=\QQ[a,b,c,d,e]$ and $\cA$ a $3\times 5$ matrix with entries
general linear forms in $R$, we got that $\Tor^R_{j}(M,M)$ and
$\Tor^R_{j}(S_{2}M,B_t)$ have different graded Betti numbers for $j=1$ and for
$j = 2$. Thus $\Tor^R_{j}(S_{2}M,B_t) \not\cong \Tor^R_{j}(M,M)$ for $j=1$ and for
$j = 2$. Note $\depth_JB_t = \dim B_t = 2$.

(4) Finally taking $R=\QQ[a,b,c,d]$ and $\cA$ as above, we got that
$\Hom(S_{-1}M,B_t)$ and $M$ 
have different graded Betti numbers. As
$\Tor^R_{3}(S_3M,B_t) \cong \Hom(S_{-1}M,B_t)$ and
$\Tor^R_{3}(S_2M,M)\cong \Hom(B_t,M) \cong M$ by \eqref{torext} below,
the hypothesis $\depth_JB_t \ge 2$ 
cannot in general be dropped.
\end{remark}

The conjecture may even hold outside $[-1, c]$, 
and inside $[-1, c]$, we may use results from \cite{KM2015}  about isomorphisms between corresponding $\Ext^i$-groups to partially prove Conjecture \ref{newconj3}. Indeed, as by Proposition \ref{resol}(ii), the $R$-dual of the $R$-free resolution of $S_r(M)$ and the $R$-free resolution of $S_{c-1-r}(M)$ coincide up to twist, we get
  \begin{equation} \label{torext} \Ext^{c-j}_R(S_{c-1-r}M,S_{s}M) \cong
\Tor^R_{j}(S_{r}M,S_{s}M)(\ell)\, ,\ \ \ \ {\rm for} \  -1 \le r, s \le c
\end{equation}
with $\ell :=\sum _{j=1}^{t+c-1}a_j-\sum _{i=1}^tb_i$, cf. Remark 3.22(ii)  (in its latest version on the arXiv of \cite{KM2015}). Moreover we proved in
Corollary 3.17 (in the latest version in arXiv) that for $s' \ge
-1$, $r'-1 \le s' \le
c$, 
\begin{equation} \label{extextn}
  \Ext^i_R(S_{r'}M,S_{s'}M) \cong  \Ext^i_R(B_t,S_{s'-r'}M) \ \ \ {\rm for\ every \ }
   i \le  \min \{a,r'- s'+c\}
  \end{equation}
%
  where $a\ge 0$ is an integer satisfying $\depth _JB_t\ge 2a+2$ (resp.
  $\depth _JB_t\ge 2a+1$ if $s'-r' \in \{-1,0,c-2\}$ and $a > 1$) with
  $J=I_{t-1}(\varphi)$. Let's use \eqref{extextn} twice, both to
  $(r',s') = (c-1-r,s)$ and to $(r',s') = (c-1-r+q,s+q)$, noticing that
  $s'-r'=r+s+1- c$ in both cases. Also $r'-1 \le s'$ yields $c-2 \le r+s$ and
  $r'- s'+c$ equals $2c-1-r-s$ in both cases. So by \eqref{extextn} we get
  \begin{equation}\label{newnew} \Ext^{i}_R(S_{c-1-r}M,S_{s}M) \cong \Ext^{i}_R(S_{c-1-r+q}M,S_{s+q}M) \ \
  {\rm for\ every \ } i \le \min \{a,2c-1-r-s\}\end{equation}
  under the assumption $c-2 \le r+s$ and $ -1 \le s, s+q \le c$. 
   Using \eqref{torext}  with $c-j$ for $i$, (\ref{newnew}) shows that \eqref{conj53eq}
 of Conjecture \ref{newconj3} holds when $\depth _JB_t\ge 2$ and 
\begin{equation} \label{extcond1}
   j \ge \max \{r+s-c+1, c+1-\frac{\depth _JB_t}2 \} \quad  {\rm and} \quad
   r+s \ge c-2.
\end{equation}
   The ``resp.-assumptions'' above give slightly more when
   $r+s+1-c \in \{-1,0,c-2\}$ and $\depth _JB_t\ge 5$, namely that
\begin{equation} \label{extcond2}
 j \ge \max \{r+s-c+1, c+\frac{1-\depth _JB_t}2 \} 
\end{equation}
imply \eqref{conj53eq}; note $\depth _JB_t\ge 5$ is included  to ensure
 $a > 1$, and $r+s\ge c-2$ is automatic as $r+s+1-c\in \{-1,0,c-2\}$. For example if $(c,r,s)=(3,1,1)$ and $\depth _JB_t=c+2$, we get \eqref{conj53eq}
for $ j\ge  1$ by \eqref{extcond2} and $-2 \le q \le 2$ from $-1\le s,s+q\le c$ which imply what we
separately prove in the first part of Lemme \ref{lem3}(ii).

Note that Conjecture \ref{newconj3} holds for
 $j \ge c$. Indeed $S_iM$ for $i \in \{s, s+q\}$ is maximally Cohen-Macaulay,
 thus has projective dimension $c$. So both $\Tor_j$-groups in
 \eqref{conj53eq} vanish for $j > c$. For $j = c$ these $\Tor_j$-groups are
 isomorphic by \eqref{torext} as $\depth _JB_t\ge 2$ yields
 $\depth _J\Hom_R(S_kM,S_iM) \ge 2$ and thus
 $\Hom(S_kM,S_{i}M)\cong \Hl _*^0(U,\widetilde{\Hom(S_kM,S_{i}M)})$,
 $U:=\Proj B_t\setminus V(J)$ for any $k \ge -1$. As $\widetilde{ S_kM}$ is
 invertible in $U$,
 $\widetilde{\Hom(S_{r'}M,S_{s}M)} \cong \widetilde{\Hom(S_{r'+q}M,S_{s+q}M)}$
 over $U$. Letting  $r' = c-1-r$ and using \eqref{torext} we get that
 \eqref{conj53eq} holds.

 Also notice, as explained in Remark \ref{IBt} that for generic determinantal
 rings we have used Macaulay2 over the field $\mathbb Q$ to get
 $\Tor^R_{1}(S_2M,B_t) \cong \Tor^R_{1}(M,M)$ for
  $2 \le c \le 5$ and $t = 3$ and also for $2 \le c \le 9$ and $t = 2$, giving
  further evidence
 to an important case of \eqref{conj53eq} of Conjecture \ref{newconj3}.

 An interesting consequence of the conjecture is that if
 $\max\{0,\frac{c+3-\depth_JB_t}{2},r+s-c\} < j \le r+s+1$ and
 $\depth _JB_t\ge 2$, then each $\Tor_j$-group of the conjecture has $\depth$
 at least $2$ at $J$. So if $U=\Proj B_t \setminus V(J)$,
 then 
   \begin{equation} \label{depthtor} \Tor^R_{j}(S_{r}M,S_{s}M)\cong \Hl
_*^0(U,\widetilde{\Tor^R_{j}(S_{r}M,S_{s}M)}).\end{equation} This is because
the conjecture implies $\Tor^R_{j}(S_{r}M,S_{s}M) \cong
\Tor^R_{j}(S_{j}M,S_{s+r-j}M) $ where the latter $\Tor_j$-group is the $j$-th
homology group of
the $R$-free resolution of $S_{j}M$ in \eqref{splice}:
$$ \to S_{0}(F)\otimes \bigwedge^tF\stackrel{\epsilon _{j}}{\longrightarrow}
\bigwedge^{j} G^* \otimes S _{0}(F^*) \rightarrow \bigwedge ^{j-1} G^* \otimes
S_1(F^*)\rightarrow \ldots \rightarrow \bigwedge^0G^*\otimes S_j(F^*) \to
S_{j}M \rightarrow 0 , $$ tensored by $N:=S_{s+r-j}M$. Since ${\epsilon _{j}}
\otimes id_{N}=0$ it follows that $ \Tor^R_{j}(S_{j}M,S_{s+r-j}M)$ is
the kernel of  $ G^* \otimes S _{0}(F^*)\otimes N \rightarrow \bigwedge ^{j-1}
G^* \otimes S_1(F^*)\otimes N $
and since $N=
S_{s+r-j}M$ is maximally Cohen-Macaulay, due to $r+s-c \le j \le r+s+1$, we
get $\depth _J N \ge 2$ and hence \eqref{depthtor}.

It is also of interest to ask if \eqref{conj53eq} holds when $c = 1$. In fact
if $\depth_JB_t \ge 2$ the answer is yes. Indeed using that
$ \Tor^R_{i}(N_1,N_2) \cong \Tor^R_{i}(N_2,N_1)$ where $N_j$, $j = 1,2$ are
finitely generated $R$-modules, there's only one non-trivial case to consider,
namely to show
\begin{equation} \label{c1n} \Tor^R_{1}(B_t,B_t) \cong \Tor^R_{1}(M,S_{-1}M)\,
  .
 \end{equation}
 Let's show \eqref{c1n} for any $c \ge 1$ assuming $\depth_JB_t \ge c+1$, and
 that all entries of $\cA$ are general of positive degree if $c \ge 2$. Firstly
 $ \Tor^R_{1}(B_t,B_t) \cong B_t \otimes I_t \cong I_t/I_t^2$ has
 $\depth_J \ge 2$ by \cite[Proposition 4.1]{KM2015} if $c \ge 2$ (and $B_t$ is
 a complete intersection if $c = 1$). Next
 $\depth_J \Tor^R_{1}(M,S_{-1}M) \ge 2$ by Remark \ref{continuation}. Thus
 \eqref{depthtor} applies to both $\Tor_1$-groups of \eqref{c1n}, and since,
 locally on $U$,
 $$\widetilde{\Tor_{1}(M,S_{-1}M)} \cong \widetilde{S_{-1}M\otimes
   \Tor^R_{1}(M,B_t)}\cong \widetilde{S_{-1}M\otimes M \otimes I_t}\cong
 \widetilde{I_t/I_t^2}$$
 by arguing as in the proof of Lemma \ref{S2M} with $S_{-1}M$ in place of $M$,
 both $\Tor_1$-groups are isomorphic, as desired, (cf.  \cite[Remarks 3.22(v)]{KM2015} of the
 latest version on the arXiv).

\begin{remark} \label{equivconj3n} \rm 
Note that using \eqref{torext} twice, we get that Conjecture~\ref{newconj3}  is equivalent to:

  Let $i < \min(c,\frac{c-3+\depth_JB_t}{2}), r, s$ and $q$ be
  integers 
  and assume $\depth_JB_t \ge 2$ and $c \ge 2$. It holds
  \begin{equation} \label{equivconj} \Ext^i_R(S_rM,S_sM) \cong \
\Ext^i_R(S_{r+q}M,S_{s+q}M)  \ \ {\rm for} \ \ r, s, r+q,s+q \in [-1, c]\, ,
\end{equation}
see \cite[Proposition 3.9]{KM2015} and its corrections (Corollary 3.17 and Remarks 3.18 and 3.22) appearing 
in the latest arXiv version of \cite{KM2015}. In
Remark 3.18, we have $\depth _J B_t=c+2$ and
we say ``we expect Corollary 3.17 to be true for $i \le c-1$ and $X$
general'' (with $X:=\Proj B_t$), but Remark \ref{tonewconj3}(1) and (2) of this paper provide examples with $c-1
= 2$ where
$\Tor^R_{1}(S_{2}M,B_t) \not\cong \Tor^R_{1}(M,M)$. Thus by \eqref{torext},
\begin{equation} \label{excor317}
\Ext^i_R(M,M) \not\cong \ \Ext^i_R(B_t,B_t)
\end{equation}
for $i = 2$, and Corollary 3.17 (cf. \eqref{extextn} with $r' = s' = 1$) is in
general not true for $i=c-1$. In Remark \ref{tonewconj3}(3) we even have an
example where also the corresponding $\Tor_2$-groups are non-isomorphic. Thus
\eqref{torext} imply that \eqref{excor317} holds also for $i = 1$ (and
$i = 2$) showing that Corollary 3.17 does not extend to cover the case
$i = c-2$ either. In all these examples, however, we have $\depth _JB_t=\dim B_t<c+2$ and what we in  Remark 3.18 say has to be interpreted as the case of a general determinantal $B_t$ of dimension $\ge c+2$ (otherwise $\depth _JB_t=c+2$ is impossible).

The equivalent version of Conjecture \ref{newconj3} stated in \eqref{equivconj} carefully clarifies all this. Indeed, if $\depth _JB_t=c+2$ (in which case $\dim B_t\ge c+2$) the assumption on the $i$ in \eqref{equivconj} is just: $i<\min(c,(c-3+\depth_JB_t)/2)=\min(c,c-1/2)$ or, equivalently $i\le c-1$ while the case $\depth_JB_t<c+2$ only ensure \eqref{equivconj} to hold for some $i<c-1$; in fact all ``counterexamples'' from Remark \ref{tonewconj3}  above to extend the conclusion in Corollary 3.17 do not satisfy $i<\min(c,(c-3+\depth_JB_t)/2)$ of \eqref{equivconj}. 
Moreover, we remark that the weak depth
condition $\depth_JB_t \ge 2$ is assumed in Conjecture \ref{newconj3} also
because we often use $\depth_JB_t \ge 2$ when we partially prove
\eqref{conj53eq} of the conjecture, or equivalently \eqref{equivconj}.

 Finally if we in \eqref{extcond1} and \eqref{extcond2} replace $j$ by $c-i$
 and $r$ by $c-1-r$ we get conditions 
 under which \eqref{equivconj} holds.
\end{remark}


\end{document}